\documentclass[10pt,journal,draftcls,onecolumn]{IEEEtran}

\usepackage{graphicx,epstopdf}
\usepackage{amsmath,amsthm,amssymb}
\usepackage{mathrsfs}
\usepackage{amssymb}
\usepackage{url}
\usepackage{color}
\usepackage{subeqnarray}
\usepackage{bbm}
\usepackage{algorithm}
\usepackage{algorithmic}
\usepackage{cite}

\newtheorem{definition}{\bf Definition}

\newtheorem{lemma}{\bf Lemma}
\newtheorem{theorem}{\bf Theorem}

\newtheorem{example}{\bf Example}
\newtheorem{remark}{\bf Remark}
\newtheorem{proposition}{\bf Proposition}
\newtheorem{corollary}{\bf Corollary}

\newcommand{\hN}{\mathcal{N}}

\newcommand{\hC}{\mathcal{C}}
\newcommand{\hS}{\mathbb{S}}
\newcommand{\hO}{\mathbb{O}}

\newcommand{\hX}{\mathbb{X}}

\newcommand{\Ex}{{\bf \mathbb{E}}}

\newcommand{\prob}{P}
\newcommand{\pr}{\mathbb{P}}
\newcommand{\conv}{{\mathrm{conv}}}
\newcommand{\dist}{\mathrm{dist}}
\newcommand{\subopt}{\mathrm{SubOpt}}

\usepackage{comment}
\usepackage{ifthen}
\newboolean{showcomments}
\setboolean{showcomments}{true}
\newcommand{\ex}[1]{  \ifthenelse{\boolean{showcomments}}{ {\bf E}#1} {}  }
\newcommand{\id}[1]{  \ifthenelse{\boolean{showcomments}}{ {\bf 1}_{#1} } {}  }
\newcommand{\eqdef}{:=}

\usepackage{ifthen}

\newcommand{\lgan}[1]{  \ifthenelse{\boolean{showcomments}}
{ \textcolor{red}{(Lingwen says:  #1)}} {}  }
\newcommand{\adam}[1]{\ifthenelse{\boolean{showcomments}}
{ \textcolor{Green}{(Adam says:  #1)}}{}}
\newcommand{\steven}[1]{\ifthenelse{\boolean{showcomments}}
{ \textcolor{Green}{(Steven says: #1)}}{}}
\newcommand{\niangjun}[1]{\ifthenelse{\boolean{showcomments}}
{ \textcolor{Green}{(Niangjun says:  #1)}}{}}
\newcommand{\addcite}[0]{\ifthenelse{\boolean{showcomments}}
{ \textcolor{blue}{(addcite)}}{}}
\newcommand{\addcites}[0]{\ifthenelse{\boolean{showcomments}}
{ \textcolor{blue}{(addcite(s))}}{}}
\newcommand{\addref}[0]{\ifthenelse{\boolean{showcomments}}
{ \textcolor{blue}{(addref)}}{}}
\newcommand{\todo}[1]{\ifthenelse{\boolean{showcomments}}
{ \textcolor{red}{(To do: #1)}} {} }

\begin{document}

\title{\Large\bf Distributed Load Balancing with Nonconvex Constraints: A Randomized Algorithm with Application to Electric Vehicle Charging Scheduling}
\author{Lingwen Gan, Ufuk Topcu, and Steven Low
}
\maketitle

\begin{abstract}
With substantial potential to reduce green house gas emission and reliance on fossil fuel, electric vehicles (EVs) have lead to a booming industry, whose growth is expected to continue for the next few decades. However, EVs present themselves as large loads to the power grid. If not coordinated wisely, the charging of EVs will overload power distribution circuits and dramatically increase power supply cost. To address this challenge, significant amount of effort has been devoted in the literature to schedule the charging of EVs in a power grid friendly way. Nonetheless, the majority of the literature assumes that EVs can be charged intermittently at any power level below certain rating, while in practice, it is preferable to charge an EV consecutively at a pre-determined power to prolong the battery lifespan. This practical EV charging constraint is nonconvex and complicates scheduling.

To schedule a large number of EVs with the presence of practical nonconvex charging constraints, a distributed and randomized algorithm is proposed in this paper. The algorithm assumes the availability of a coordinator which can communicate with all EVs. In each iteration of the algorithm, the coordinator receives tentative charging profiles from the EVs and computes a broadcast control signal. After receiving this broadcast control signal, each EV generates a probability distribution over its admissible charging profiles, and samples from the distribution to update its tentative charging profile. 

We prove that the algorithm converges almost surely to a charging profile in finite iterations. The final charging profile (that the algorithm converges to) is random, i.e., it depends on the realization. We characterize the final charging profile---a charging profile can be a realization of the final charging profile if and only if it is a Nash equilibrium of the game in which each EV seeks to minimize the inner product of its own charging profile and the aggregate electricity demand. Furthermore, we provide a uniform suboptimality upper bound, that scales $O(1/n)$ in the number $n$ of EVs, for all realizations of the final charging profile.
\end{abstract}

\section{Introduction}
Electric vehicles (EVs) are propelled by electric motors powered by rechargeable battery packs, while conventional vehicles run on internal combustion engines. EVs enjoy a higher energy efficiency, more environmental friendliness, a superior performance (in terms of noise, acceleration, maintenance, etc.), and less energy dependency on foreign fossil fuels than conventional vehicles \cite{doe2014all}. Therefore EVs are promoted by the US government with tax incentives---a \$7500 federal income tax credit for every EV (small neighborhood EVs excluded) purchased in or after 2010 \cite{doe2014federal}. Besides government support, EVs are much more fuel economic, e.g., the 2014 Honda Fit EV is estimated to cost \$0.87 electricity per 25 miles drive \cite{doe2014compare}, while the 2014 Honda Fit is estimated to cost \$3.25 gas per 25 miles drive.\footnote{Assume that a 2014 Honda Fit burns 1 gallon gas per 30 miles drive, and each gas costs \$3.90.} Such government support and fuel economy have lead to a booming EV industry: over 170000 highway-capable EVs have been sold in the US from 2008 to 2013, and 16 EV models rom 9 major manufacturers are available in the US market by March 2014 \cite{wiki2014plug}. This trend is forecast to speed up as major vehicle manufacturers consecutively announce their EV plans \cite{BMW, Nissan, Tesla}.

A multitude of challenges, not only economical but also technical, have to be overcome to welcome the new era of EVs. On the economical side, EVs are currently sold at a much higher price than conventional vehicles, mainly due to the expense of battery packs. This is starting to change since battery prices are dropping by 20\%--30\% each year \cite{Hans2013mindboggling}, and therefore EVs are likely to have a competitive price in the near future. Technical challenges are more fundamental: among other things, EVs have a limited range due to small battery capacities, and it takes hours to get the batteries fully charged. Except for Tesla Model S, all other currently available EV models have a range of at most slightly exceeding 100 miles \cite{doe2014compare}. Typically, it takes 3.6--12 hours to fully charge a battery, depending on the battery capacity and the charging technology \cite{doe2014compare}.

The technical challenge that we seek to solve in this paper is the integration of EVs to the power grid. While EVs present themselves as large loads to the power delivery circuits, the power delivery circuits are designed something like five decades ago without bearing in mind the prosperous of EVs. Furthermore, EVs are likely to be charged after their owners arrive home from work since large-scale public charging facilities are still to be constructed. If not coordinated wisely, the charging of EVs may lead to coincidence peaks in electricity demand \cite{kelly2009analyzing}. Consequently, power transmission/distribution lines carry much larger currents and power transformers are loaded much more heavily. As a result, circuit device lifespan will be greatly reduced \cite{roe2009power}, and voltages throughout the network may suffer from severe deviations from nominal values \cite{clement2009coordinated}---note that electricity appliances like air conditioners and televisions only work under well-regulated voltages.

On the other hand, if the charging of EVs is well coordinated, not only will some of the integration challenges be mitigated, but also ancillary services can be provided to stabilize the voltages and frequency \cite{quinn2010effect,lopes2009smart,gan2013real}, and power supply cost can be reduced \cite{verzijlbergh2012network}. For example, if the charging of EVs is coordinated such that the aggregate electricity demand is flat, then high peak electricity demand can be avoided. Consequently, circuit devices will be least loaded and their lifespans will be enlarged; voltages throughout the network will be less fluctuating and easier to regulate; electricity demand is easier to follow since it ramps up and down much slower, and therefore expensive power plants for following the demand can be reduced, bringing down the power supply cost. As another example, the charging of EVs can be coordinated to respond to network contingencies such as the unexpected shutting down of a power plant. In such cases, the charging of some EVs can be postponed to rebalance the electricity supply and demand, which is essential for preventing electricity blackouts.

\subsection*{Literature review}
A variety of algorithms have been proposed in the literature to coordinate the charging of EVs, ranging from centralized algorithms where a coordinator makes decisions for the EVs \cite{qian2011modeling, verzijlbergh2012network, sundstrom2010planning, richardson2012optimal} to distributed algorithms where each EV makes their own decisions \cite{ma2013decentralized, gan2013optimal, caramanis2009management, rotering2011optimal}. Note that distributed algorithms may still require coordinators to facilitate the communication among EVs.

Centralized algorithms are mainly for cost-benefit analysis purposes when the number of EVs is large, since the computation burden would be too heavy for a single computation unit. In \cite{verzijlbergh2012network}, a large number of operational distribution networks in The Netherlands are investigated to quantify the impact of EV charging on various network levels. Results show that controlled charging can reduce infrastructure update investment by half over uncontrolled charging. In \cite{qian2011modeling}, uncontrolled charging and smart charging of EVs are compared empirically to highlight that a second peak electricity demand during night can be avoided by adopting smart charging. In \cite{sundstrom2010planning} and \cite{richardson2012optimal}, centralized linear programmings are proposed to compute the optimal charging profiles of EVs. In \cite{sundstrom2010planning}, the linear programming aims to minimize the power supply cost subject to circuit capacity constraints and vehicle owner's requirements. The linear programming in \cite{richardson2012optimal} further includes a power network physical model that captures voltage deviations.

By distributing the computation burden among different units, a distributed algorithm is more suitable for scenarios where a large number of EVs need to be coordinated. In \cite{ma2013decentralized}, a distributed algorithm is proposed to schedule the charging of EVs such that the aggregate electricity demand is made flat. It is proved in \cite{ma2013decentralized} that when the EVs are identical, the obtained aggregate electricity demand will be optimal in the sense that it is as flat as possible. This notion of optimal aggregate electricity demand is formalized in \cite{gan2013optimal}. Besides, a different distributed algorithm is proposed in \cite{gan2013optimal} to schedule the charging of EVs. Furthermore, the algorithm proposed in \cite{gan2013optimal} always obtains optimal aggregate electricity demand. The algorithms proposed in \cite{caramanis2009management, rotering2011optimal} take another perspective: instead of trying to flatten the aggregate electricity demand, they seek to minimize the charging cost given a pre-determined electricity price profile by solving a dynamic programming at each EV.
 
All aforementioned works make the assumption that the battery of an EV can be charged intermittently at any power level below certain rating, e.g., if Single-Phase Level II charging is deployed, then the rating is 3.3kW \cite{ipakchi2009grid} and an EV can be charged at any power level from 0 to 3.3kW when it is plugged in. However in practice, to prolong the lifespan of batteries, it is preferably to charge the batteries consecutively at a fixed pattern until they get fully charged \cite{james2003electric}. Then, the only flexibility in scheduling the charging of EVs is to determine when the charging of each EV starts. Such flexibility can no longer be described by linear constraints. Moreover, one can show that such flexibility imposes nonconvex constraints on the scheduling of EV charging, calling for a different approach.

\subsection*{Summary of contributions}
This paper proposes a distributed algorithm to schedule the charging of EVs in order to shape the aggregate electricity demand (e.g., to flatten the aggregate electricity demand or to track some target demand profile), in the presence of convex or nonconvex charging constraints (recall that charging constraints are nonconvex if the only flexibility in scheduling is the time when an EV starts charging). In particular, our contributions are threefold.

First, in the special case where charging constraints are convex, we improve the convergence rate of the distributed algorithm proposed in \cite{gan2013optimal} by updating the charging profiles of EVs at different speeds based on their total energy consumption. A load balancing problem that aims to minimize the $\ell_2$ norm of the aggregate load profile is formulated in Section \ref{subsec: lb}, and the optimal EV charging scheduling problem formulated in Section \ref{subsec: app to ev} is an example of the load balancing problem. In the special case where charging constraints are convex, a distributed algorithm is presented in Section \ref{subsec: special case} to compute the optimal charging schedules of EVs.

This algorithm, referred to as Algorithm \ref{algorithm: convex}, is similar to the gradient projection algorithm (GPA) proposed in \cite{gan2013optimal}. Both algorithms are iterative, assume the availability of a coordinator that communicates with all EVs, and share the same communication pattern. In each iteration of Algorithm GPA, the coordinator receives tentative charging profiles from EVs and computes the partial derivatives of the objective function with respect to each tentative charging profile. The partial derivatives are identical since the objective function only depends on the aggregate of charging profiles, and is therefore symmetric (the interchange of two arbitrary charging profiles does not change the objective value). Then, the coordinator broadcasts the common partial derivative to all EVs. After receiving the common partial derivative, each EV updates its tentative charging profile by first taking a step along the negative direction of the partial derivative and then projecting the resulting charging profile back to the set of admissible charging profiles. It has been proved in \cite{gan2013optimal} that the sequence of tentative charging profiles obtained by Algorithm GPA always converges to the set of optimal charging profiles.

Algorithm \ref{algorithm: convex} improves over Algorithm GPA in that the step sizes, taken by the EVs to update their tentative charging profiles, are chosen to speed up the convergence rate. The intuition of choosing the step sizes is that, EVs with larger total energy requests should update their charging profiles faster since they have a larger set of admissible charging profiles to optimize over. In particular, the step sizes are chosen in proportional to the total energy requests of EVs.

Second, in the general case where charging constraints are nonconvex but certain technical conditions hold, we propose a randomized distributed algorithm for scheduling the charging of EVs. This algorithm, referred to as Algorithm \ref{algorithm: nonconvex}, defers from Algorithm \ref{algorithm: convex} only in how EVs update their tentative charging profiles. Recall that in Algorithm \ref{algorithm: convex}, EVs update their tentative charging profiles by first taking a step in the direction of decreasing the objective value, and then projecting the resulting charging profile back to the set of admissible charging profiles. However when charging constraints are nonconvex, the set of admissible charging profiles can be difficult to project onto. Moreover, one can prove that any deterministic algorithm is guaranteed to diverge in the simplest setup where EVs are identical, unless the communication pattern is different from Algorithm \ref{algorithm: convex}. This highlights the necessity of a randomized algorithm to handle the nonconvex charging constraints.

In Algorithm \ref{algorithm: nonconvex}, to update its tentative charging profile, an EV first computes a charging profile $\hat{x}$ in the convex hull of admissible charging profiles, then finds a probability distribution over admissible charging profiles that has expectation $\hat{x}$, and finally samples from the probability distribution to update its tentative charging profile. Noteworthy, Algorithms \ref{algorithm: convex} and \ref{algorithm: nonconvex} share the same coordinator operations: the coordinator receives tentative charging profiles and computes a broadcast signal---the common partial derivative of the objective function with respect to individual charging profiles. Such compatibility between Algorithms \ref{algorithm: convex} and \ref{algorithm: nonconvex} allows for the joint scheduling of EVs with convex and nonconvex constraints: for EVs with convex constraints, run the EV algorithm in Algorithm \ref{algorithm: convex}; for EVs with nonconvex constraints, run the EV algorithm in Algorithm \ref{algorithm: nonconvex}.

The derivation of Algorithm \ref{algorithm: nonconvex} is based on the martingale theory \cite{grimmett2001probability}. Loosely speaking, a supermartingale is a stochastic process that decreases in expectation, and all lower bounded supermartingales converge almost surely. Algorithm \ref{algorithm: nonconvex} is designed such that the objective values obtained in consecutive iterations form a supermartingale. More specifically, the randomness in Algorithm \ref{algorithm: nonconvex} comes from the step in which each EV samples from a probability distribution to update its tentative charging profile. The probability distributions are designed such that on average, the objective value decreases from the previous iteration.

Third, we provide performance guarantees for Algorithm \ref{algorithm: nonconvex}. A charging profile $x$ of all EVs is called stationary if the sequence of tentative charging profiles obtained by Algorithm \ref{algorithm: nonconvex} stops updating once the sequence reaches $x$. We prove that stationary charging profiles exist and provide a characterization for them. In particular, a charging profile is stationary if and only if it is a Nash equilibrium of the game in which each EV seeks to minimize the inner product of its own charging profile and the aggregate electricity demand. We also provide a uniform upper bound on the suboptimality of stationary charging profiles. The upper bound scales $O(1/n)$ as the number $n$ of EVs increases. Hence, when the number $n$ of EVs is large, all stationary charging profiles are nearly optimal.

We prove that the sequence of tentative charging profiles obtained by Algorithm \ref{algorithm: nonconvex} converges to stationary charging profiles. More specifically, we prove that the sequence of tentative charging profiles converges to the set of stationary charging profiles. Furthermore, in cases where each EV only has finitely many admissible charging profiles, the sequence of tentative charging profiles converges to a single but random stationary charging profile, i.e., the sequence may converge to different stationary charging profiles in different realizations. In numerical simulations, Algorithm \ref{algorithm: nonconvex} converges fast to nearly optimal charging profiles. In particular, after 10 iterations, the charging profile obtained by Algorithm \ref{algorithm: nonconvex} incurs a suboptimality less than 3\%, for a wide range of EV penetration levels.

\section{Problem formulation} \label{sec: problem formulation}
We study the offline charging scheduling of EVs in this paper. In particular, we assume that EVs are available for negotiation at the beginning of a time horizon, on their charging profiles over the horizon. The goal is to compute an EV charging schedule that flattens the aggregate electricity demand while respecting the constraints of EVs.

The EV charging scheduling problem can be cast in a more general setting, where there are workloads to be distributed among several workers. The goal is to distribute the workloads such that each worker has more or less the same amount of workload. In the context of real applications, workloads can be software processes, Internet traffic, data to be transmitted, etc., and workers can be cpu cores, routers, communication channels, etc..

In this section, we first formulate a load balancing problem in Section \ref{subsec: lb}, and then describes its application to EV charging scheduling in Section \ref{subsec: app to ev}.

\subsection{The load balancing problem}
\label{subsec: lb}
Consider the scenario where a number of loads need to be served over a time horizon. The loads can be divided into two categories: elastic loads whose service time and/or rate can be adjusted, and inelastic loads whose service time and rate are fixed and given. The load balancing problem considered in this paper seeks to schedule the service time and rate of elastic loads such that about the same aggregate service rate is required at each time.

Let $T$ denote the length of the time horizon, and assume that the time horizon starts at 0 without loss of generality. In practice, $T$ could be a day. The service rates of inelastic loads can be aggregated into a single process $b(t)$ where $t\in[0,T]$. Since inelastic loads can only be served at fixed time and rate, the process $b(t)$ is not adjustable. Hence, the process $b(t)$ is called base load in the literature.

Let $n$ denote the total number of elastic loads and index them by $1,2,\ldots,n$. For each elastic load $i$, let $x_i(\cdot):[0,T]\rightarrow\mathbb{R}$ denote its service rate profile, i.e., a service of rate $x_i(t)$ is required by elastic load $i$ at time $t$. The service time/rate of an elastic load is adjustable, and this can be captured by a region where the service rate profile of the elastic load can be optimized over. More specifically, for each elastic load $i$, there exists a region $\hX_i$ such that the flexibility in $x_i$ is captured by
	\[ x_i \in \hX_i. \]
Every element in $\hX_i$ is called an admissible service rate profile. Assume that each admissible service rate profile requires the same amount of total service, i.e., for each elastic load $i$, there exists $X_i\in\mathbb{R}$ such that for any admissible service rate profile $y_i\in\hX_i$, one has
	\[ \int_0^T y_i(t)dt = X_i. \]
The quantity $X_i$ is called the total service requirement.

%

The goal is to flatten the aggregate service rate profile
	\[ d:=b+\sum_{i=1}^n x_i. \]
This is captured by minimizing the variance
	\[ V(d)\eqdef\frac{1}{T}\int_0^T \left(d(t)-\frac{1}{T}\int_0^T d(\tau)d\tau \right)^2dt =\frac{1}{T}\int_0^T d^2(t)dt - \mu_d^2 \]
of $d$, where the quantity
	\[ \mu_d := \frac{1}{T}\left(\int_0^Tb(t)dt + \sum_{i=1}^n X_i\right) \]
is the time average of $d$. Note that $\mu_d$ is independent of $x_i$, and therefore minimizing the variance $V(d)$ of $d$ over $x:=(x_1,x_2,\ldots,x_n)$ is equivalent to minimizing the squared $\ell_2$ norm
	\[ \|d\|_2^2 = \int_0^T d^2(t)dt \]
of $d$ over $x$.

To summarize, the load balancing problem is formulated as follows.
	\begin{subequations}
	\label{LB}
    	\begin{align}
    	\textbf{LB: }\min~~ & \int_0^T \left(b(t) + \sum_{i=1}^nx_i(t) \right)^2dt \label{LB: obj} \\
    	\text{s.t.}~~ & x_i\in \mathbb{X}_i,\quad i=1,\ldots,n. \label{nonconvex}
    	\end{align}
	\end{subequations}
Note that the problem LB can be nonconvex since the sets $\mathbb{X}_i$ can be nonconvex.

The following assumption is made on the constraint sets $\hX_i$ throughout this work. To state the assumption, let $\mathcal{C}[0,T]$ denote the set of real-valued continuous functions defined on $[0,T]$, and equip $\mathcal{C}[0,T]$ with the $\ell_\infty$ norm $\|\cdot\|_\infty:\mathcal{C}[0,T]\mapsto\mathbb{R}_+$ defined as
	\[ \|f\|_\infty := \max_{t\in[0,T]}\{ |f(t)| \}, \qquad \forall f\in\hC[0,T]. \]
We assume that the set $\hX_i$ is compact throughout this paper.
	
\begin{remark}
The algorithms and analysis developed in this paper apply to more general objective functions. In particular, the objective function in \eqref{LB: obj} can be generalized to
	\begin{equation}
	\label{LB: obj generalize}
	\int_0^T C_t\left(b(t) + \sum_{i=1}^n x_i(t)\right)dt
	\tag{\ref{LB: obj}'}
	\end{equation}
where the function $C_t:\mathbb{R}\mapsto\mathbb{R}$ represents the cost at time $t$ and has a uniform upper bound $\beta_t$ on its second derivatives $C_t''(\tau)$ over $\tau\in[0,T]$ for $t\in[0,T]$, i.e., there exists $\beta_t\in\mathbb{R}$ for $t\in[0,T]$ such that
	\[ C_t''(\tau) \leq \beta_t, \qquad \forall t,\tau\in[0,T]. \]
Note that if $C_t(x)=x^2$ for all $t\in[0,T]$, then \eqref{LB: obj generalize} simplifies to \eqref{LB: obj}. Also note that the objective of shaping the aggregate service profile $d$ to track a target profile $G:[0,T]\mapsto\mathbb{R}$ can be captured by setting
	\[ C_t(x) = (x-G(t))^2 \]
for $t\in[0,T]$ in \eqref{LB: obj generalize}.
\end{remark}

\subsection{Application to EV charging}
\label{subsec: app to ev}
The EV charging scheduling problem can be considered as an application of the problem LB. In a power network, there are both electricity demand due to EVs and electricity demand due to other loads such as air conditioners, televisions, and lights. The charging of EVs can be postponed as long as EVs get fully charged by certain deadlines, and therefore EVs can be considered as elastic loads. Control of other loads including air conditioners, televisions, and lights is more evolved and outside the scope of this paper. Hence, we model those loads as inelastic.

At the beginning of a time horizon, the electricity utility negotiates with the EVs on their charging profiles over the time horizon, in order to achieve certain system level objectives, e.g., to flatten the aggregate electricity demand over the time horizon or to shape the aggregate electricity demand to track certain target profile. In this paper, we focus on the objective of flattening the aggregate electricity demand, and the results generalize straightforwardly to shaping the aggregate electricity demand to track certain target profile.

The objective of flattening the aggregate electricity demand is chosen for the following reasons. A flat electricity demand has the lowest peak and therefore induces the least burden to circuit devices. Consequently, the lifespan of circuit devices will be prolonged. Besides, the voltages in the network fluctuate as electricity demand varies. Hence, a flat electricity demand causes the least amount of fluctuation in network voltages and simplifies voltage regulation. At last, electricity generation and demand need to be balanced in almost real time since electricity storage is limited. In practice, such balance is achieved by controlling the electricity generation of some power plants called ``spinning reserves'' to follow the demand. These power plants need to ramp up and down fast enough to follow the fluctuation in electricity demand, and are usually expensive. If electricity demand is flat, then fewer ``spinning reserves'' will be necessary and electricity supply cost can be significantly reduced.

Adopt the notations introduced in Section \ref{subsec: lb}, i.e., let $[0,T]$ denote the time horizon, let $b(t)$ denote the aggregate electricity demand of inelastic loads, index the EVs by $1,2,\ldots,n$, and let $x_i(t)$ denote the charging profile of EV $i$ for $i=1,2,\ldots,n$.

	\begin{figure}[!htbp]
      	\centering
      	\includegraphics[scale=0.3]{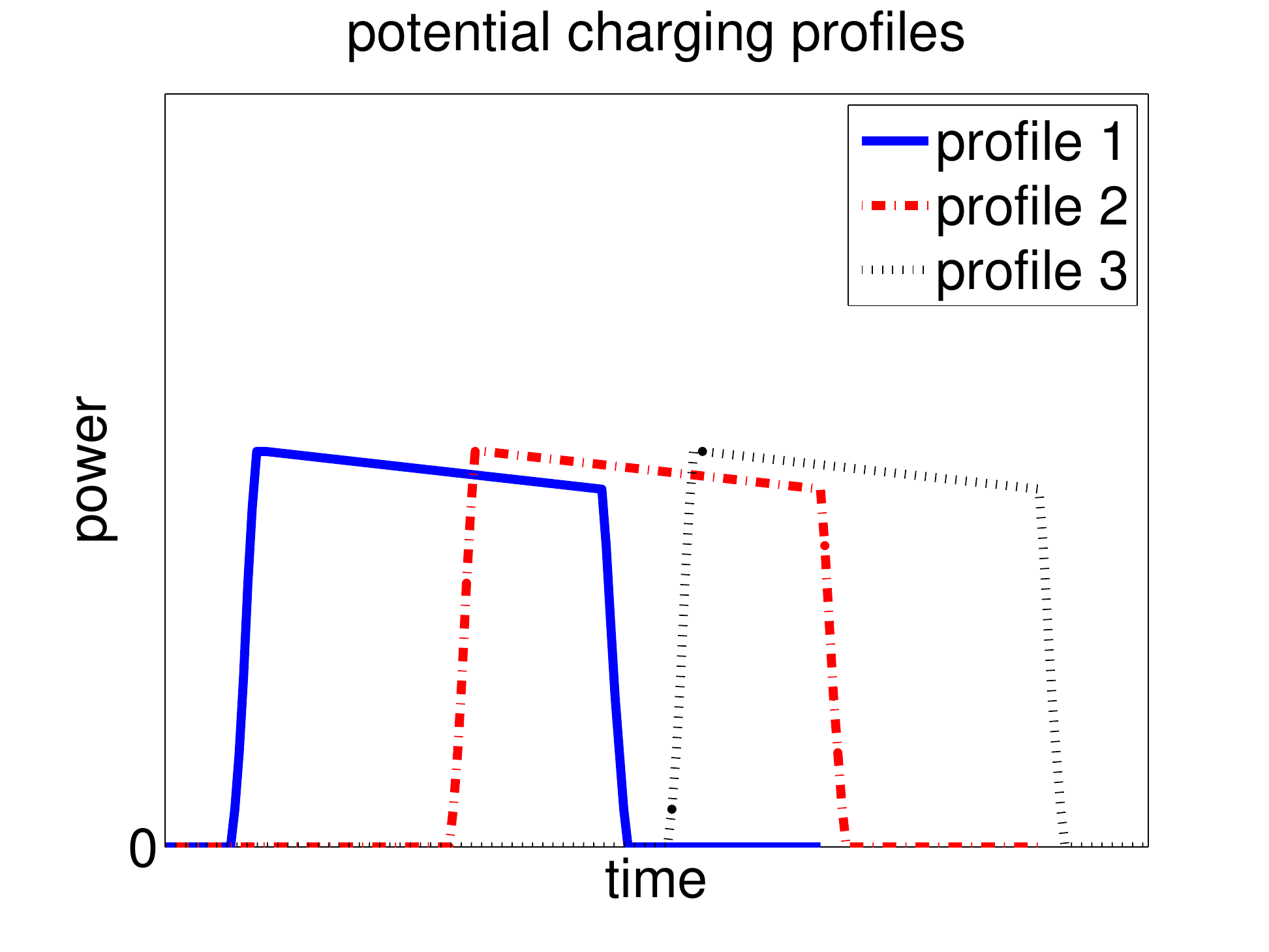}
      	\caption{Hypothetical charging profiles of an EV. They are time-shifted versions of each other.}
      	\label{fig: charge}
	\end{figure}

Charging an EV at certain fixed pattern can prolong the lifespan of its battery \cite{james2003electric}. If such patterns are enforced in the scheduling problem, then admissible charging profiles $y_i\in\hX_i$ of an EV $i$ are time-shifted versions of each other as illustrated in Fig. \ref{fig: charge}. In this case, the set $\hX_i$ of admissible charging profiles is not convex and therefore the problem LB is difficult to solve, but elements $y_i$ in $\hX_i$ satisfy the following properties:
	\begin{itemize}
	\addtolength{\itemindent}{0.4cm}
	\item[A1)] Charging rate is uniformly bounded, i.e., there exists $\overline{x}_i\in\mathbb{R}$ such that
		\[ |y_i(t)| \leq \overline{x}_i, \qquad \forall t\in[0,T], ~ \forall y_i\in\hX_i. \]
	\item[A2)] Ramp rate is uniformly bounded, i.e., there exists $L_i\in\mathbb{R}$ such that
		\[ |y_i(s)-y_i(t)| \leq L_i|s-t|, \qquad \forall s,t\in[0,T], ~ \forall y_i\in\hX_i. \]
	\item[A3)] Total energy request is fixed, i.e., there exists $X_i\in\mathbb{R}$ such that
		\[ \int_0^Ty_i(t)dt=X_i, \qquad \forall y_i\in \hX_i. \]
	\item[A4)] The $\ell_2$ norm of charging profile is fixed, i.e., there exists $Y_i\in\mathbb{R}$ such that
		\[ \int_0^Ty_i^2(t)dt=Y_i, \qquad \forall y_i\in \mathbb{X}_i. \]
	\end{itemize}

\begin{remark}
\label{remark: compact}
Let $\conv(\hX_i)$ denote the convex hull of $\hX_i$, i.e.,
	\[ \conv(\hX_i) = \left\{ \sum_{k=1}^m \theta_ky_k \mid m\in\mathbb{Z}^+, ~\theta_k\geq0\text{ for all }k, 
	~\sum_{k=1}^m \theta_k=1, ~y_k\in\hX_i \text{ for all }k \right\}, \]
then every element $y_i\in\conv(\hX_i)$ also satisfies A1, A2, and A3. Therefore, the set $\conv(\hX_i)$ is uniformly bounded (by A1) and equi-continuous (by A2). The set $\conv(\hX_i)$ is also closed since we assume $\hX_i$ to be compact throughout this paper. Therefore, the set $\conv(\hX_i)$ is compact in the metric space $(\mathcal{C}([0,T]), \|\cdot\|_\infty)$ by the Arzela-Ascoli Theorem \cite[Theorem 11.18]{carothers2000real}.
\end{remark}

\begin{remark}
\label{remark: simplex}
The sets $\hX_i$'s can be simplified to be finite (a compact set can be approximated by a set of finite elements to arbitrary precision) in order to ease the computation of the algorithm developed in Section \ref{subsec: algorithm nonconvex}. More specifically, we restrict the times when an EV starts charging to finite choices. For instance, if an EV can start charging at any time from 20:00 to 3:00, then we enforce the EV to start charging only at 20:00, 20:05, $\ldots$, 3:00. Consequently, the set $\hX_i$ of admissible charging profiles is finite.

Let $m$ denote the number of elements in $\hX_i$, and denote these elements by $y_1,y_2,\ldots, y_{m}$. Every element $z\in\conv(\hX_i)$ is a convex combination of $y_1,y_2,\ldots, y_{m}$, i.e., there exists $\theta=(\theta_1,\theta_2,\ldots,\theta_{m})$ such that
	\[ \theta_k\geq0\text{ for all }k, ~\sum_{k=1}^{m}\theta_k=1, ~z=\sum_{k=1}^{m}\theta_ky_k. \]
It is straightforward to show that the map
	\[ f: z \mapsto \theta \]
is continuous and one-to-one from $\conv(\hX_i)$ to the probability simplex
	\[ \mathbb{P}^m = \left\{ p=(p_1,p_2,\ldots,p_m) \mid p_i\geq0 \text{ for all }i, ~ \sum_{i=1}^m p_i = 1 \right\}, \]
i.e., $\conv(\hX_i)$ is homeomorphic to $\mathbb{P}^m$.
\end{remark}

\section{Distributed algorithms} \label{sec: algorithm}
The focus of this paper is the design and analysis of a distributed algorithm to solve the problem LB. In the special case where constraint sets $\hX_i$ are convex, we propose a distributed algorithm Algorithm \ref{algorithm: convex} in Section \ref{subsec: special case}. In each iteration of Algorithm \ref{algorithm: convex}, a coordinator broadcasts a control signal in order to guide the loads to update their tentative service rate profiles. In the general case where constraint sets $\hX_i$ are not necessarily convex but their convex hulls $\conv(\hX_i)$ are, we propose a randomized distributed algorithm Algorithm \ref{algorithm: nonconvex} in Section \ref{subsec: algorithm nonconvex}. Algorithm \ref{algorithm: nonconvex} shares the same coordinator operations as Algorithm \ref{algorithm: convex}, but differs in the load operations.

\subsection{Special case: convex constraint sets $\hX_i$}
\label{subsec: special case}
When $\hX_i$ is convex for $i\in\hN$, the problem LB can be solved by a gradient projection approach as in our prior work \cite{gan2013optimal}. Here, we modify the algorithm proposed in \cite{gan2013optimal} by allowing loads to update their tentative service rate profiles at different speeds in order to increase the convergence rate.

\vspace{-0.15in}
\subsection*{\bf Infrastructure requirement}
The algorithm assumes the availability of a coordinator that can broadcast control signals to all loads and receive feedbacks from the loads. In the application to EV charging, the coordinator can be an electricity utility or a lower level aggregator, and the communication between the EVs and the coordinator can be enabled by for example the Home Area Networks.

\vspace{-0.15in}
\subsection*{\bf Algorithm statement}
To state the algorithm, define the inner product $\left\langle \cdot,\cdot\right\rangle: \hC[0,T]\times\hC[0,T]\rightarrow \mathbb{R}$ for two continuous functions $f$, $g$ on $[0,T]$ as
	\[ \left\langle f,g \right\rangle := \int_0^T f(t)g(t)dt, \]
and define the $\ell_2$ norm $\|\cdot\|:\hC[0,T]\rightarrow\mathbb{R}$ for a continuous function $f$ on $[0,T]$ as
	\[ \|f\| := \sqrt{\langle f,f \rangle}. \]
	
\begin{algorithm}[!h]
\caption{Constraint set $\hX_i$ is compact and convex for $i\in\hN$}
\label{algorithm: convex}
\begin{algorithmic}
    \REQUIRE Horizon length $T$, base load $b$, constraint sets $\hX_i$ for $i\in\hN$, componentwise strictly positive parameter $c=(c_1,\ldots,c_n)\in\mathbb{R}^n$, tolerance $\epsilon>0$.

    \ENSURE A service profile $x=(x_1,\ldots,x_n)$.

    \STATE {\bf Initialization}\\
    $\quad$ Each load $i\in\hN$ initializes $x_i^{(0)} \leftarrow 0$;\\
    $\quad$ $k \leftarrow 1$;
    \STATE {\bf Repeat}\\
    $\quad$ The coordinator collects $x_i^{(k-1)}$ from all loads and computes a control signal
    	\begin{equation}\label{price update}
	g^{(k)} \leftarrow \frac{b+\sum_{i=1}^nx_i^{(k-1)}} {\sum_{i=1}^nc_i};
	\end{equation}
    $\quad$ Each load $i\in\hN$ receives the control signal $g^{(k)}$, and updates its service rate profile as
    	\begin{equation}\label{rate update}
	 x_i^{(k)}\leftarrow\underset{x_i\in\hX_i}{\mathrm{argmin}}~ 2c_i\left\langle g^{(k)}, x_i \right\rangle + \left\| x_i-x_i^{(k-1)} \right\|^2;
	\end{equation}
    $\quad$ $k\leftarrow k+1$;
    \STATE {\bf Until} $k>2$ and $\left\|g^{(k-1)}-g^{(k-2)}\right\|<\epsilon$.
    \STATE{\bf Return} $x \leftarrow x^{(k-1)}$.
\end{algorithmic}
\end{algorithm}

The algorithm is iterative and works as follows. In each iteration $k=1,2,\ldots$:
\begin{itemize}
\item The coordinator receives tentative service profiles $x_i^{(k-1)}$ from the loads, and computes the normalized aggregate service profile $g^{(k)}$ (normalize by $\sum_{i=1}^n c_i$). Then, the coordinator broadcasts $g^{(k)}$ as a control signal to all loads in order to coordinate the loads in updating their tentative service profiles.

\item Each load updates its tentative service profile by solving \eqref{rate update} after receiving the broadcast control signal $g^{(k)}$ from the coordinator. The objective in \eqref{rate update} is the sum of two terms $\left\langle g^{(k)}, x_i \right\rangle$ and $\left\| x_i-x_i^{(k-1)} \right\|^2$ weighted by $2c_i$ and 1 respectively. The first term $\left\langle g^{(k)}, x_i \right\rangle$, if interpreting $g^{(k)}$ as the time-varying service price profile, is the service cost associated with service profile $x_i$. The second term $\left\| x_i-x_i^{(k-1)} \right\|^2$ penalizes the deviation from $x_i$ to the service profile $x_i^{(k-1)}$ computed in the previous iteration $k-1$. The second term is added in \eqref{rate update} to ensure the convergence of Algorithm \ref{algorithm: convex}.
\end{itemize}

The parameters $c_i$ weigh the two terms in \eqref{rate update} and can be chosen appropriately to speed up the convergence rate of Algorithm \ref{algorithm: convex}. In particular, assume that A3 holds, i.e., the total service $X_i$ required by a load $i$ is a constant that does not depend on the choice of admissible service profile $y_i\in\hX_i$ for $i\in\hN$, and $c_i$ is set as $c_i=X_i$ for $i\in\hN$. The reason will become clear after Theorem \ref{thm: rate convex}.

\vspace{-0.15in}
\subsection*{\bf Optimality}
Algorithm \ref{algorithm: convex} solves the problem LB in the special case where constraint sets $\hX_i$ are convex and compact. To state the result, for two service rate profiles $x=(x_1,\ldots,x_n)$ and $x'=(x_1',\ldots,x_n')$, define the distance $d(x,x')$ between $x$ and $x'$ as
	\[ d(x,x'):= \sqrt{\sum_{i=1}^n \left\|x_i-x_i'\right\|^2}. \]
Let $\hO^*$ denote the set of optimal solutions of LB, and define the distance $d(x,\hO^*)$ from a service rate profile $x$ to the set $\hO^*$ as
	\[ \dist(x,\hO^*):=\inf_{x^*\in\hO^*} d(x,x^*). \]
If $\dist(x,\hO^*)=0$, then $x$ is in the closure $\overline{\hO^*}$ of $\hO^*$. If $\dist(x,\hO^*)$ is small, then $x$ is ``close'' to $\hO^*$.

\begin{theorem}\label{thm: convergence convex}
If $\hX_i$ is compact and convex for $i\in\hN$, then the sequence $\{x^{(0)}, x^{(1)}, \ldots, x^{(k)}, \ldots\}$ of service profiles computed in Algorithm \ref{algorithm: convex} converges to the set $\hO^*$ of optimal solutions of LB as $k\rightarrow \infty$, i.e.,
	\[ \dist(x^{(k)},\hO^*)\rightarrow 0 \text{ as } k\rightarrow \infty. \]
\end{theorem}

Theorem \ref{thm: convergence convex} implies that Algorithm \ref{algorithm: convex} obtains an optimal service profile. It is proved in Appendix \ref{app: convergence convex}. Note that no further assumption is made on the parameter $c$ in Theorem \ref{algorithm: convex}, i.e., Algorithm \ref{algorithm: convex} eventually solves the problem LB for any parameter $c$ that is componentwise strictly positive. However, a good choice of $c$ can significantly improve the convergence rate of Algorithm \ref{algorithm: convex}.

\vspace{-0.15in}
\subsection*{\bf Convergence rate}
The parameter $c$ is picked for a faster convergence of Algorithm \ref{algorithm: convex}, and its choice is based on the following observations. First, $c_i$ affects how fast a load $i$ updates its tentative service profile $x_i$: with a larger $c_i$, more weight is placed on the first summand in \eqref{rate update} and therefore the second summand in \eqref{rate update}, a penalty term on service profile updates, is less impactful. Consequently, load $i$ tends to make larger updates in $x_i$.

Second, the second summand in \eqref{rate update} scales quadratically with total service $X_i$, while $\left\langle g^{(k)}, x_i \right\rangle$ in the first summand scales linearly with total service $X_i$. Hence, if $c_i$ is set proportional to $X_i$, then the two summands in \eqref{rate update} are scale invariant, i.e., the two summands in \eqref{rate update} remain the same ratio as $X_i$ varies over a wide range. According to the first observation, this choice of $c_i$ enforces loads with larger total service $X_i$ to take larger steps in their service profile updating, such that the tentative service profiles of different loads converge to the final service profiles at a ``uniform'' rate.

The following pathological example explains the benefit of setting $c_i$ proportional to $X_i$.
\begin{example}
\label{example: 3 loads}
Consider two instances I1 and I2 of the problem LB. In Instance I1, there are 3 identical loads 0, 1, 2, each with constraint set $\hX$ and total service $X$. It is not difficult to verify that if we set $c_0=c_1=c_2$, then Algorithm \ref{algorithm: convex} obtains an optimal service profile within 1 iteration.

In Instance I2, there are only two loads 1' and 2'. While load 2' is the same as load 2, i.e., load 2' has constraint set $\hX$ and requests total service $X$, load 1' is effectively the aggregate of load 0 and load 1, i.e., load 1' has constraint set $\hX_{1'} = 2\hX$ and requests total service $X_{1'} = 2X$. If we still set $c_{1'} = c_{2'}$, then Algorithm \ref{algorithm: convex} no longer converges within 1 iteration. This is because load 1' has a larger feasible set $\hX_{1'}$ to explore than load 2', but can only update its tentative service profiles $x_{1'}$ at the same speed as load 2'. To remedy this, one can set $c_{1'}=2c_{2'}$. Interestingly, with this choice of $c_{i'}$, one can prove that $x_{1'}^{(k)}=2x_0^{(k)}=2x_1^{(k)}$ and $x_{2'}^{(k)}=x_2^{(k)}$ for $k=0,1,2,\ldots$, i.e., Algorithm \ref{algorithm: convex} behaves equivalently in Instances I1 and I2 in the sense that load $1'$ is effectively the aggregate of load 0 and load 1. Hence, Algorithm \ref{algorithm: convex} still converges within 1 iteration, highlighting the benefit of setting $c_i$ proportional to $X_i$.
\end{example}

We can generalize Example \ref{example: 3 loads} to multiple-load setups as stated in the following theorem.
\begin{theorem}
\label{thm: rate convex}
Assume that $\hX_i$ is compact and convex, and that $c_i=X_i>0$ for each $i$. Consider the following two instances of the problem LB.
\begin{itemize}

\item[I1)] There are $n+1$ loads $0,1,\ldots,n$, and load 0 and load 1 are identical.

\item[I2)] There are $n$ loads $1',2',\ldots,n'$ where load 1' is the aggregate of load 0 and load 1 while loads $2'\sim n'$ are the same as loads $2\sim n$, i.e., 
	\[ \hX_{1'} = 2\hX_0 = 2\hX_1, \quad \hX_{i'} = \hX_i \text{ for }i=2,3,\ldots,n. \]

\end{itemize}
The sequence $\left\{ [x_i^{(0)}]_{i=0}^n,     [x_i^{(1)}]_{i=0}^n,  \ldots,     [x_i^{(k)}]_{i=0}^n,   \ldots \right\}$ of tentative service profiles computed in Algorithm \ref{algorithm: convex} in Instance I1 is equivalent to the sequence $\left\{ [x_{i'}^{(0)}]_{i=1}^n,     [x_{i'}^{(1)}]_{i=1}^n,  \ldots,     [x_{i'}^{(k)}]_{i=1}^n,   \ldots \right\}$ of tentative service profiles computed in Algorithm \ref{algorithm: convex} in Instance I2 in the sense that
	\begin{equation}\label{split}
	x_{1'}^{(k)} = 2x_0^{(k)} = 2x_1^{(k)}, \qquad x_{i'}^{(k)}=x_i^{(k)} \text{ for }i=2,3,\ldots,n
	\end{equation}
for $k=0,1,2,\ldots$.
\end{theorem}
Theorem \ref{thm: rate convex} implies that the aggregation of identical loads into a bigger load does not change the convergence rate of Algorithm \ref{algorithm: convex} provided that $c_i=X_i$ for all $i$. The theorem is proved in Appendix \ref{app: rate convex}. Hence, we set $c_i=X_i$ for all $i$ in Algorithm \ref{algorithm: convex}.

\vspace{-0.15in}
\subsection*{\bf Connection to prior work}
Algorithm \ref{algorithm: convex} generalizes the distributed algorithm proposed in \cite{gan2013optimal} by introducing a parameter $c$ to regulate the tentative service profile updates of different loads. In particular, if $c_i=1$ for all $i$, then Algorithm \ref{algorithm: convex} reduces to the algorithm proposed in \cite{gan2013optimal}. The freedom of choosing $c$ in Algorithm \ref{algorithm: convex} can be exploited to improve the convergence rate as exemplified in Example \ref{example: 3 loads}. In this paper, we set $c$ as $c_i=X_i$ for all $i$.

\subsection{General case: nonconvex constraint sets $\hX_i$}
\label{subsec: algorithm nonconvex}
The problem LB can be solved by Algorithm \ref{algorithm: convex} in the special case where constraint sets $\hX_i$ are convex. However, the constraint sets $\hX_i$ are not convex for many applications including the application to EV charging, and therefore Algorithm \ref{algorithm: convex} does not directly apply. In this section, we seek to modify Algorithm \ref{algorithm: convex} such that certain types of nonconvex constraint sets $\hX_i$ can be dealt with.

Motivated by the application to EV charging, we assume that each $\hX_i$ is compact and satisfies A1, A2, A3, and A4 in the rest of this section. It follows from Remark \ref{remark: compact} that $\conv(\hX_i)$ is compact for each $i$.

\vspace{-0.15in}
\subsection*{\bf Infrastructure requirement}
We aim at an algorithm with the same infrastructure requirement and information flow pattern as Algorithm \ref{algorithm: convex} so that loads with convex and nonconvex constraint sets can be scheduled simultaneously. In particular, we assume the availability of a coordinator that is capable of broadcasting control signals to all loads and receiving feedbacks.

\begin{figure}[!htbp]
      \centering
      \includegraphics[scale=0.5]{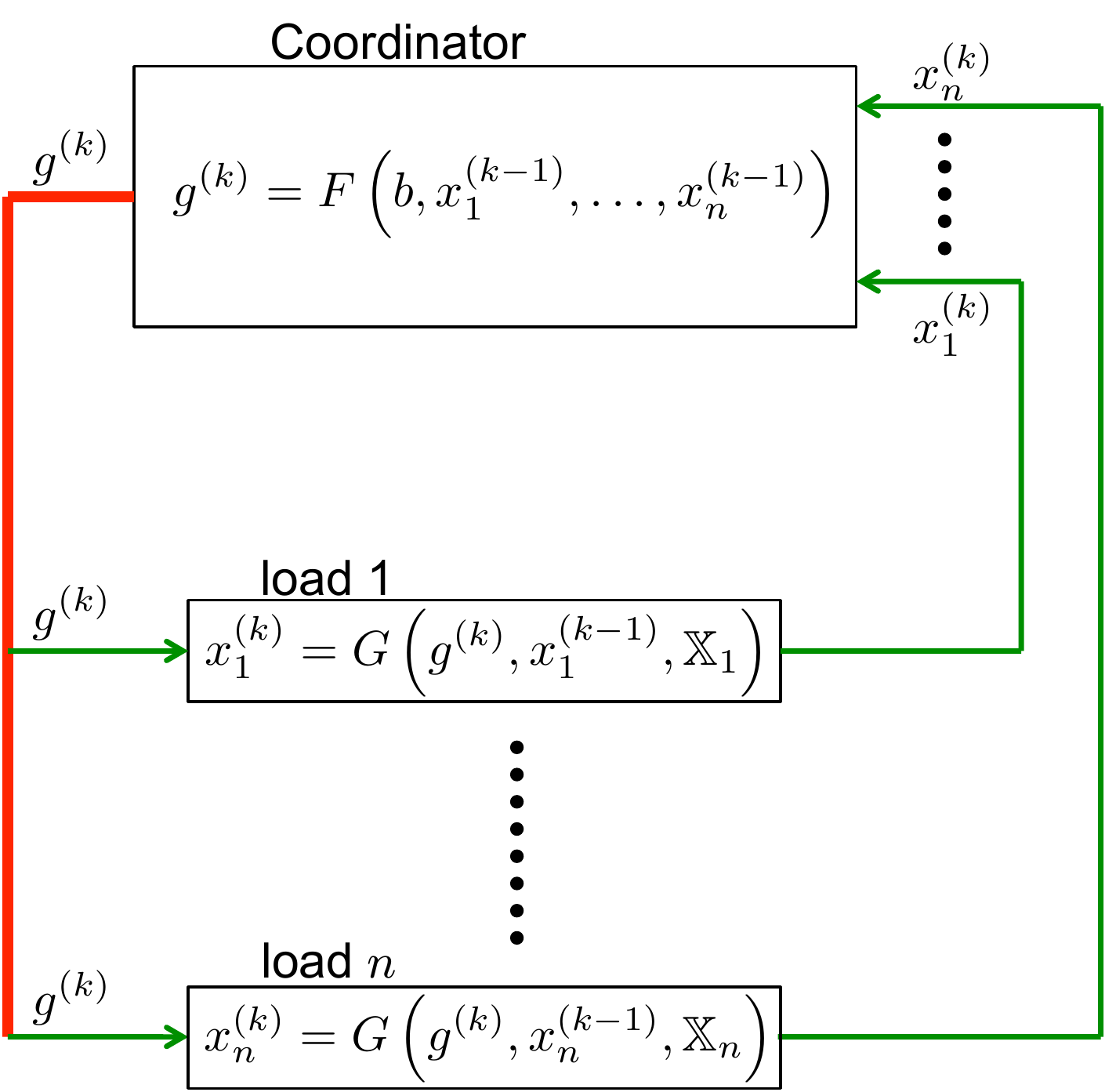}
      \caption{Information flow pattern in the proposed iterative, distributed decision-making process. The coordinator knows the base load $b$ and each load $i$ knows its admissible service profiles $\hX_i$.}
      \label{fig:pattern}
\end{figure}

The information flow pattern of Algorithm \ref{algorithm: convex} (as well as the algorithm to be designed) is described in Fig. \ref{fig:pattern} where a coordinator and multiple loads exchange information in multiple iterations to agree on a service schedule $x=(x_1,x_2,\ldots,x_n)$. In each iteration $k$, the coordinator computes $F$ to obtain a control signal $g^{(k)}$ based on the base load profile $b$ and tentative service profiles $x_1^{(k-1)},\ldots,x_n^{(k-1)}$ computed by loads in the previous iteration, and each load $i$ computes $G$ to update its tentative service profile $x_i^{(k)}$ based on the control signal $g^{(k)}$, the tentative service profile $x_i^{(k-1)}$ it computes in the previous iteration, and its set of admissible service profiles $\hX_i$. Here $F$ and $G$ are the functions to be designed.

In order to schedule the loads with convex and nonconvex constraints simultaneously, the algorithm to be designed must share the same coordinator operations, i.e., the same function $F$, as Algorithm \ref{algorithm: convex}. However, the load operations, i.e., the function $G$, can be different for the algorithm to be designed and Algorithm \ref{algorithm: convex}. 

When there are both loads with convex constraints and loads with nonconvex constraints in the system, the coordinator does not need to behave differently and simply computes $F$ in each iteration, while each load, depending on whether it has convex constraints or nonconvex constraints, computes the $G$ for Algorithm \ref{algorithm: convex} or the $G$ for the algorithm to be designed, respectively.

\vspace{-0.15in}
\subsection*{\bf Failure of deterministic algorithms}
Any algorithm with the information flow pattern in Fig. \ref{fig:pattern} and a deterministic function $G$ either diverges or converges to synchronized solutions in the setting where all loads are identical.

\begin{proposition}[Proposition 1, \cite{gan2012stochastic}]
\label{pro:deter}
If $\hX_1=\hX_2=\ldots=\hX_n$ and $G$ is a deterministic function, then any algorithm with the information flow pattern in Fig. \ref{fig:pattern} obtains synchronized service profiles in all iterations, i.e.,
	\[ x_1^{(k)}=x_2^{(k)} = \cdots = x_n^{(k)}, \qquad k\geq0. \]
\end{proposition}

Proposition \ref{pro:deter} implies that if loads are identical, then deterministic algorithms obtain synchronized service profiles. When $\hX_i$'s are convex, synchronized service profiles can still be optimal and that's why Algorithm \ref{algorithm: convex}, a deterministic algorithm, can be optimal in that setting. However when $\hX_i$'s are nonconvex, even the best synchronized service profile can be (and are usually) far from optimal. This highlights the necessity of randomized algorithms.

\begin{remark}
Proposition \ref{pro:deter} only discusses deterministic algorithms using the information flow pattern in Fig. \ref{fig:pattern}. It is possible to find a deterministic algorithm that flattens the aggregate load using other information flow patterns. For instance, consider the following information flow pattern. In each iteration $k$, the coordinator uses the most recently calculated $x_1,\ldots,x_n$ to compute the broadcast signal $g$, and only one load (in turn) updates its tentative service profile. With this information flow pattern, there exists deterministic algorithms that flatten the aggregate load, however, at a much higher communication overhead and a much longer delay.
\end{remark}

\vspace{-0.15in}
\subsection*{\bf Idea of randomized algorithms}

A randomized algorithm adopts a randomized function $G$ and works as follows. In each iteration $k$, each load $i$ first computes a probability distribution $P_i^{(k)}$ over its admissible service profiles $\hX_i$ and then samples from the probability distribution $P_i^{(k)}$ to update its tentative service profile $x_i^{(k)}$. The algorithmic challenge is the design and computation of a probability distribution $P_i^{(k)}$ that ensures the convergence of a randomized algorithm.

The probability distributions $P_i^{(k)}$ are designed such that the objective values in consecutive iterations of the algorithm decrease in expectation. With this property and lower boundedness of the objective values, there is a prominent theorem saying that the objective values almost surely converge as iterations continue.

\vspace{-0.15in}
\subsection*{\bf Algorithm statement}
To state the algorithm, let $\Theta(A)$ denote the set of probability distributions over a set $A$. For example, if $A$ has $m$ elements, then any probability distribution $P$ over $A$ can be described by a point $p=(p_1,p_2,\ldots,p_m)$ in the probability simplex $\mathbb{P}^m$, i.e., there exists a one-to-one mapping between $\Theta(A)$ and $\mathbb{P}^m$.
Given $f:A\rightarrow \mathbb{R}$ and a probability distribution $P$ over $A$, let
	\[ \Ex_P[f]:=\int_A f(y)d\prob \]
denote the expectation of $f$ with respect to the probability distribution $\prob$.

\begin{algorithm}[!h]
\caption{Constraint set $\hX_i$ is compact and satisfies A1, A2, A3, A4 for $i\in\hN$}
\label{algorithm: nonconvex}
\begin{algorithmic}
    \REQUIRE Horizon length $T$, base load $b$, constraint sets $\hX_i$ for $i\in\hN$, componentwise strictly positive parameter $c=(c_1,\ldots,c_n)\in\mathbb{R}^n$, tolerance $\epsilon>0$.
    \ENSURE A service profile $x=(x_1,\ldots,x_n)$.
    \STATE {\bf Initialization}\\
    $\quad$ Each load $i\in\hN$ initializes $x_i^{(0)} \leftarrow 0$;\\
    $\quad$ $k \leftarrow 1$;
    \STATE {\bf Repeat}\\
    $\quad$ The coordinator collects $x_i^{(k-1)}$ from all loads and computes a control signal
    	\begin{equation}\label{price update2}
	g^{(k)} \leftarrow \frac{b+\sum_{i=1}^nx_i^{(k-1)}}{\sum_{i=1}^nc_i};
	\end{equation}
    $\quad$ Each load $i\in\hN$ receives the control signal $g^{(k)}$, computes a probability distribution $P_i^{(k)}$ over $\hX_i$ as
    	\begin{equation}\label{expectation}
	 P_i^{(k)} \leftarrow \underset{P_i\in\Theta(\hX_i)}{\mathrm{argmin}} 
	 ~2c_i\left\langle \frac{g^{(k)}\sum_{j=1}^nc_j-x_i^{(k-1)}}{\sum_{j\neq i}c_j}, ~\Ex_{P_i}[x_i]\right\rangle
	 + \left\| \Ex_{P_i}[x_i] - x_i^{(k-1)} \right\|^2,
	\end{equation}
    $\quad$ and samples from $P_i^{(k)}$ to update a tentative service profile $x_i^{(k)}\in\hX_i$;\\
    $\quad$ $k\leftarrow k+1$;
    \STATE {\bf Until} $k>2$ and $\left\|g^{(k-1)}-g^{(k-2)}\right\|<\epsilon$.
    \STATE{\bf Return} $x \leftarrow x^{(k-1)}$.
\end{algorithmic}
\end{algorithm}

Algorithm \ref{algorithm: nonconvex} is similar to Algorithm \ref{algorithm: convex}, but with different load operations. In each iteration $k$ of Algorithm \ref{algorithm: nonconvex}, a load $i$ does not directly compute an $x_i^{(k)}\in\hX_i$ as in Algorithm \ref{algorithm: convex}, but first computes a probability distribution $P_i^{(k)}$ over $\hX_i$ by solving \eqref{expectation} and then samples from $P_i^{(k)}$ to update its tentative service profile $x_i^{(k)}$. The objective in \eqref{expectation} is chosen such that the objective values evaluated in consecutive iterations of Algorithm \ref{algorithm: nonconvex} decrease in expectation.

The objective in \eqref{expectation} differs from that in \eqref{rate update} in two aspects. First, the rate profile $x_i\in\hX_i$ in \eqref{rate update} is replaced by the expectation $\Ex_{P_i}[x_i]$ of a probability distribution $P_i$ over $\hX_i$ in \eqref{expectation}. Second, the first summand in the objective of \eqref{expectation} has a negative term $-x_i^{(k-1)}$ and a multiplicative factor $\eta:=\sum_{j=1}^n c_j / \sum_{j\neq i}c_j$ in front of $g^{(k)}$ and . Note that $x_i^{(k-1)}$ is negligible in comparison with $b+\sum_{j=1}^n x_j^{(k-1)}=g^{(k)}\sum_{j=1}^nc_j$, so the negative term $-x_i^{(k-1)}/\sum_{j\neq i}c_j$ is negligible. Besides, if there is no single $c_i$ that dominates the sum $\sum_{j=1}^nc_j$, then the multiplicative factor $\sum_{j=1}^n c_j / \sum_{j\neq i}c_j$ in front of $g^{(k)}$ is approximately 1. To summarize,
	\begin{equation}
	\label{approximation}
	\frac{g^{(k)}\sum_{j=1}^nc_j-x_i^{(k-1)}}{\sum_{j\neq i}c_j} \approx g^{(k)},
	\end{equation}
i.e., the objective in \eqref{expectation} is approximately equal to the objective in \eqref{rate update} with $x_i$ substituted by $\Ex_{P_i}[x_i]$.

The computation of a probability distribution $P_i^{(k)}$ can be difficult unless a simple representation of $P_i^{(k)}$ is available. As discussed in Remark \ref{remark: simplex}, in the application to EV charging, the set $\hX_i$ can be reduced to a finite set in order to ease the computation of $P_i^{(k)}$. In this case, $P_i^{(k)}$ can be represented by its probability mass function over $\hX_i$, i.e., let $m$ denote the number of elements in $\hX_i$, then there is a one-to-one mapping between $\Theta(\hX_i)$ (where $P_i^{(k)}$ takes values in) and the probability simplex $\mathbb{P}^m$. Moreover, the map $f : \Theta(\hX_i) \mapsto \conv(\hX_i)$ defined by
	\[ f(P_i) = \Ex_{P_i}[x_i] \]
is one-to-one. Hence, computing a $P_i^{(k)}\in\Theta(\hX_i)$ is equivalent to finding an $\Ex_{P_i^{(k)}}[x_i]\in\conv(\hX_i)$.

\begin{remark}\label{rem:N}
Load $i$ needs to know $C:=\sum_{j=1}^n c_j$ in order to compute \eqref{expectation}, and the coordinator can broadcast $C$ together with $g^{(k)}$ to provide this piece of information. Alternatively, load $i$ can also use the approximation in \eqref{approximation} so that it no longer needs to know $C$. In this case, the objective load $i$ minimizes is the same as the objective in \eqref{rate update} with $x_i$ substituted by $\Ex_{P_i}[x_i]$, i.e., each load, whether it has convex or nonconvex constraint, minimizes the same objective.
\end{remark}

\vspace{-0.15in}
\subsection*{\bf Stationary service profiles and optimality}
We will show that the sequence of tentative service profiles obtained by Algorithm \ref{algorithm: nonconvex} converges to certain nearly optimal service profiles. To start with, we characterize these service profiles and show that they are nearly optimal.

\begin{definition}
\label{def: stationary}
Assume that constraint set $\hX_i$ is compact and satisfies A1, A2, A3, A4 for $i\in\hN$. A service profile $x=(x_1,\ldots,x_n)$ is stationary for Algorithm \ref{algorithm: nonconvex}, if $x^{(k-1)}=x$ implies $x^{(k)}=x$ with probability 1 for $k\geq1$.
\end{definition}
That is, stationary service profiles are those once hit, tentative service profiles stop updating. Note that stationary service profiles do not exist for many randomized algorithms including simulated annealing \cite{van1987simulated} and genetic algorithm \cite{goldberg1989genetic}. In these two widely known algorithms, the probability
	\begin{equation}
	\label{escape}
	\pr\{x^{(k)}\neq x \mid x^{(k-1)}=x\}>0
	\end{equation}
is nonzero for any $k\geq1$ and any $x$, and therefore stationary service profiles do not exist for these two algorithms.

The reason simulated annealing and genetic algorithm satisfy \eqref{escape} is that they seek to escape from local optimal points, and this objective is achieved at the cost of escaping from global optimal points as well. The methodology is that, the probability of obtaining a global optimal point increases via a so-called ``cooling''/``evolving'' process as iterations continue. Though the cooling/evolving process is usually very slow, eventually a global optimal point is obtained with a strictly positive probability.

Algorithm \ref{algorithm: nonconvex} seeks to find a nearly optimal service profile within few iterations, instead of having a slow cooling (evolving) process to find an optimal service profile (which is NP hard). At the cost of small suboptimality (Theorem \ref{thm: suboptimality}), Algorithm \ref{algorithm: nonconvex} has stationary service profiles (Theorem \ref{thm: nash}) and converges to these stationary service profiles quickly---we stop Algorithm \ref{algorithm: nonconvex} after the first 20 iterations in the simulations provided in Section \ref{sec:simulation}.

A characterization of the stationary service profiles of Algorithm \ref{algorithm: nonconvex} is provided in the following theorem.

\begin{theorem}\label{thm: nash}
Assume that constraint set $\hX_i$ is compact and satisfies A1, A2, A3, A4 for $i\in\hN$. A service profile $x=(x_1,\ldots,x_n)$ is stationary for Algorithm \ref{algorithm: nonconvex} if and only if it is a Nash-equilibrium of the game in which each load $i\in\hN$ seeks to minimize $\left\langle b+\sum_{j=1}^nx_j, x_i\right\rangle$ over $x_i\in\hX_i$.
\end{theorem}

Theorem \ref{thm: nash} implies that the stationary service profiles of Algorithm \ref{algorithm: nonconvex} happen to be the Nash-equilibria of a game. It is proved in Appendix \ref{app: nash}. In this game, each load $i\in\hN$ aims to minimize the inner product of its own service profile $x_i$ and the aggregate service profile of all loads $b+\sum_{j=1}^nx_j$. Hence at a Nash equilibrium, load $i$ chooses to be serviced at times when the aggregate load $b+\sum_{j=1}^nx_j$ is low. Therefore, stationary service profiles are likely to be close to flat. Indeed, the following theorem provides an upper bound on the suboptimality of stationary service profiles of Algorithm \ref{algorithm: nonconvex}. Before stating the theorem, recall the definition of $Y_i$ in A4.

\begin{theorem}
\label{thm: suboptimality}
Assume that constraint set $\hX_i$ is compact and satisfies A1, A2, A3, A4 for $i\in\hN$. Let $x^s$ be a stationary service profile for Algorithm \ref{algorithm: nonconvex} and $x^*$ be a global optimum of the problem LB, then
	\begin{equation}
	\label{upper bound}
	\left\|b+\sum_{i=1}^nx_i^s\right\|^2 - \left\|b+\sum_{i=1}^nx_i^*\right\|^2 \leq 4\sum_{i=1}^n Y_i.
	\end{equation}
If further, elements in $\hX_i$ are nonnegative, i.e., $y_i\geq0$ componentwise for all $i\in\hN$ and all $y_i\in\hX_i$, then
	\begin{equation}
	\label{upper bound 2}
	\left\|b+\sum_{i=1}^nx_i^s\right\|^2 - \left\|b+\sum_{i=1}^nx_i^*\right\|^2 \leq 2\sum_{i=1}^n Y_i.
	\end{equation}
\end{theorem}

Theorem \ref{thm: suboptimality} implies that the suboptimality of a stationary service profile is uniformly upper bounded by constant $4\sum_{i=1}^n Y_i$, which scales $O(n)$ as the number $n$ of loads increases. The theorem is proved in Appendix \ref{app: suboptimality}. Noting that the optimal objective value scales $O(n^2)$ as $n$ increases, the suboptimality ratio defined as
	\begin{equation}
	\label{subopt}
	\subopt(x^s) := \frac{\left\|b+\sum_{i=1}^nx_i^s\right\|^2 - \left\|b+\sum_{i=1}^nx_i^*\right\|^2}{\left\|b+\sum_{i=1}^nx_i^*\right\|^2}
	\leq \frac{ 4\sum_{i=1}^n Y_i }{\left\|b+\sum_{i=1}^nx_i^*\right\|^2}
	\end{equation}
has an upper bound that scales $O(1/n)$ as $n$ increases. Hence, all stationary service profiles of Algorithm \ref{algorithm: nonconvex} are close to optimal when $n$ is sufficiently large. When $\hX_i$'s are composed of nonnegative elements, the suboptimality ratio upper bound can be improved to
	\begin{equation}
	\label{subopt2}
	\subopt(x^s) \leq \frac{ 2\sum_{i=1}^n Y_i }{\left\|b+\sum_{i=1}^nx_i^*\right\|^2}.
	\end{equation}

\vspace{-0.15in}
\subsection*{\bf Martingale and convergence}
Algorithm \ref{algorithm: nonconvex} is designed such that the objective values obtained in consecutive iterations decrease in expectation, since then there is a martingale theory ensuring the convergence of objective values. To state the result, let 
	\[ L_k:= \left\| b + \sum_{i=1}^nx_i^{(k)} \right\|^2 \]
denote the objective value obtained in iteration $k$ of Algorithm \ref{algorithm: nonconvex} for $k\geq1$.

\begin{proposition}\label{proposition: supermartingale}
Assume that constraint set $\hX_i$ is compact and satisfies A1, A2, A3, A4 for $i\in\hN$. The stochastic process $\{L_k:k\geq1\}$ is a supermartingale with respect to the process $\{x^{(k)}:k\geq1\}$, i.e.,
	\[ \Ex \left[ L_k \mid x^{(k-1)},x^{(k-2)},\ldots,x^{(0)} \right] \leq L_{k-1}, \quad k=2,3,\ldots. \]
\end{proposition}

Proposition \ref{proposition: supermartingale} implies that $\{L_k:k\geq1\}$ is a nonnegative supermartingale. It is proved in Appendix \ref{app: supermartingale}. Martingale theory tells us that a nonnegative supermartingale almost surely converges \cite{grimmett2001probability}, as stated in the following corollary.

\begin{corollary}
Assume that constraint set $\hX_i$ is compact and satisfies A1, A2, A3, A4 for $i\in\hN$. The sequence $L_0, L_1, \ldots, L_k, \ldots$ converges almost surely as $k\rightarrow\infty$.
\end{corollary}

In fact, the objective in \eqref{expectation} is designed such that $\{L_k:k\geq1\}$ forms a nonnegative supermartingale and therefore convergences. However, the convergence of objective values does not imply the convergence of tentative service profiles, though the reverse is true. Hence, we prove the following stronger result regarding the convergence of tentative service profiles. To state the result, let $\hS$ denote the set of stationary service profiles of Algorithm \ref{algorithm: nonconvex}.

\begin{theorem}\label{thm: convergence}
Assume that constraint set $\hX_i$ is compact and satisfies A1, A2, A3, A4 for $i\in\hN$. The set $\hS$ of stationary service profiles of Algorithm \ref{algorithm: nonconvex} is nonempty, and the sequence $x^{(0)}, x^{(1)}, \ldots, x^{(k)}, \ldots$ of tentative service profiles computed in Algorithm \ref{algorithm: nonconvex} converges almost surely to the set $\hS$ as $k$ tends to infinity, i.e.,
    \[ \dist(x^{(k)},\hS) \overset{\mathrm{a.s.}}\longrightarrow 0 \text{ as } k\rightarrow \infty. \]
\end{theorem}
Theorem \ref{thm: convergence} implies that the sequence of tentative service profiles computed in Algorithm \ref{algorithm: nonconvex} converges to the set of stationary service profiles. The theorem is proved in Appendix \ref{app: convergence}. However, convergence to a set does not imply convergence to a point since there can be for example limit cycles \cite[Chapter 7.0]{strogatz2006nonlinear}. In the special case where $\hX_i$'s have finitely many elements, convergence to a point can be established as in the following Corollary.

\begin{corollary}
\label{cor: convergence to points}
Assume that constraint set $\hX_i$ is compact and satisfies A1, A2, A3, A4 for $i\in\hN$. If $\hX_i$ has finitely many elements for $i\in\hN$, then the sequence $x^{(0)}, x^{(1)}, \ldots, x^{(k)}, \ldots$ of tentative service profiles computed in Algorithm \ref{algorithm: nonconvex} converges almost surely to
a random service profile $x^\infty$ that takes values in $\hS$, i.e.,
	\[ x^{(k)} \overset{\mathrm{a.s.}}\longrightarrow x^\infty \text{ as } k\rightarrow \infty \]
where $x^{\infty}$ is a random variable that takes values in $\hS$.
\end{corollary}
Corollary \ref{cor: convergence to points} implies that the sequence of tentative service profiles computed in Algorithm \ref{algorithm: nonconvex} converges to a random stationary service profile. The corollary is proved in Appendix \ref{app: convergence to points}. In the application to EV charging, the sets $\hX_i$'s are finite, and therefore Algorithm \ref{algorithm: nonconvex} converges almost surely to a single (random) stationary service profile.

\section{Case Studies} \label{sec:simulation}
We evaluate Algorithm \ref{algorithm: nonconvex} numerically in the application to EV charging in this section. In particular, we investigate how fast Algorithm \ref{algorithm: nonconvex} converges, and how close to optimal are the charging profiles obtained by Algorithm \ref{algorithm: nonconvex} after a certain number of iterations.

\begin{figure*}[!h]
      	\centering
      	\includegraphics[scale=0.3]{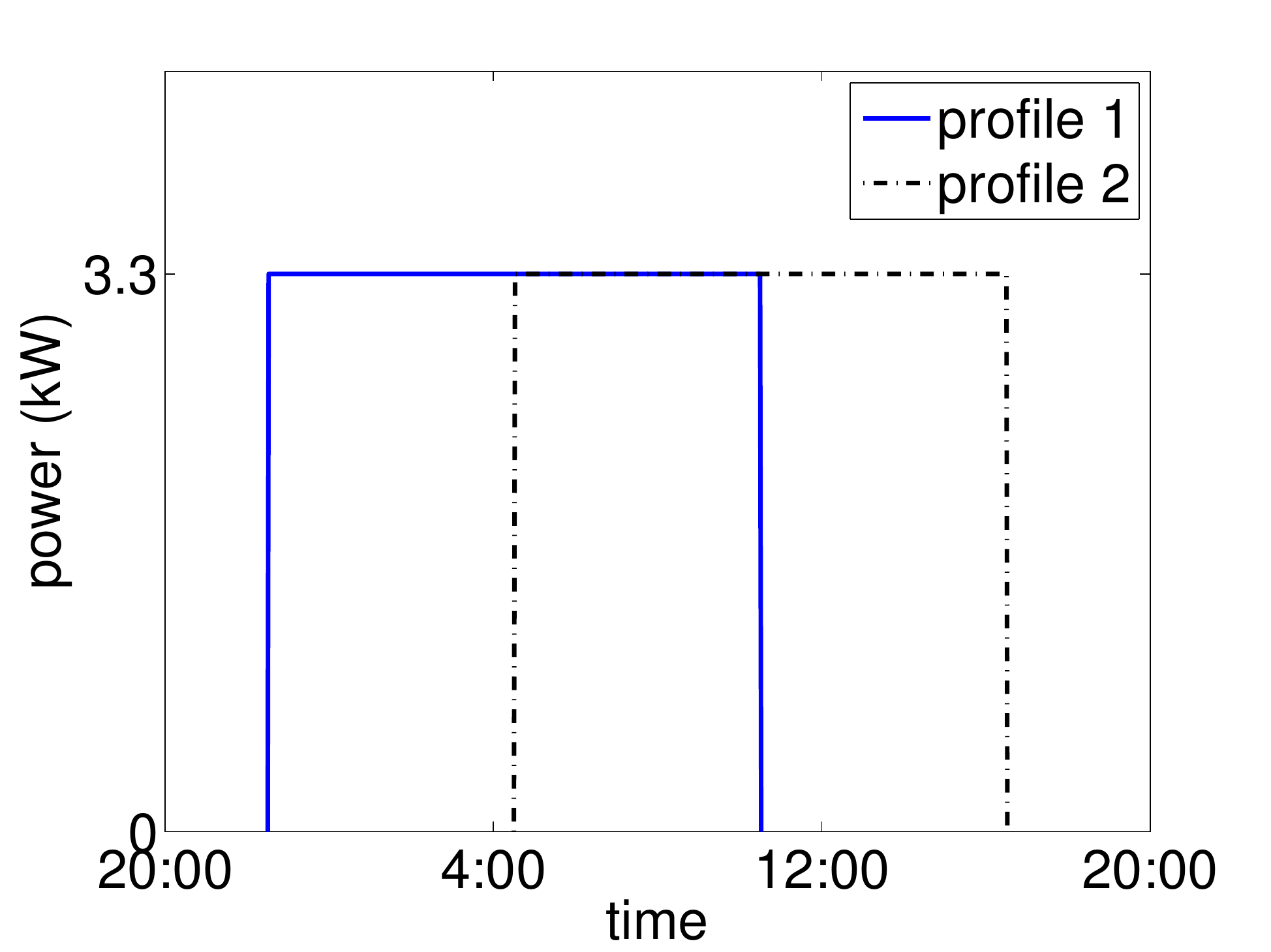}
      	\caption{Two examples of admissible charging profiles of an EV. The EV needs to charge its battery at 3.3kW power consecutively for 4 hours.}
      	\label{fig: profile}
    	\end{figure*}

The scheduling horizon $T$ is considered to be $24$ hours, and discretized into 96 slots, each of 15 minutes. The base load is set to match the average residential load profile in the service area of Southern California Edison from 20:00 on 02/13/2011 to 20:00 on 02/14/2011 \cite{SCE}. We consider different penetration levels of EVs. For simplicity, we assume that each EV is charged according to the pattern in Fig. \ref{fig: profile}---it consumes 3.3kW power consecutively for 4 hours \cite{ipakchi2009grid}. We also assume that each EV can start charging at 20:00, 20:15, $\ldots$, 16:00---as long as it can finish the 4-hour charging before the end of the time horizon.

\subsection{Convergence rate}\label{sec:simu:finite convergence}
Recall that $x^{(k)}$ denotes the tentative charging profile computed in iteration $k$ of Algorithm \ref{algorithm: nonconvex}, and note that the stochastic process $\{x^{(k)}:k\geq1\}$ is a Markov chain \cite[Chapter 6]{grimmett2001probability}. Define \emph{escape probabilities}
	\[ P_\mathrm{escape}^{(k)}:=\pr\{x^{(k)}\neq x^{(k-1)} \mid x^{(k-1)}\}, \quad k\geq1. \]
If $P_\mathrm{escape}^{(k)}=0$, then $x^{(k-1)}$ is a stationary charging profile; otherwise, it is straightforward to verify that $1/P_\mathrm{escape}^{(k)}$ is the expected number of iterations it takes before tentative charging profile gets updated. For example, if $P_\mathrm{escape}^{(k)}=0.5$, then on average a tentative charging profile update happens after 2 iterations. At such an updating speed, we may want to stop the iterations in Algorithm \ref{algorithm: nonconvex}.

\begin{figure*}[!h]
      \centering
      \includegraphics[scale=0.5]{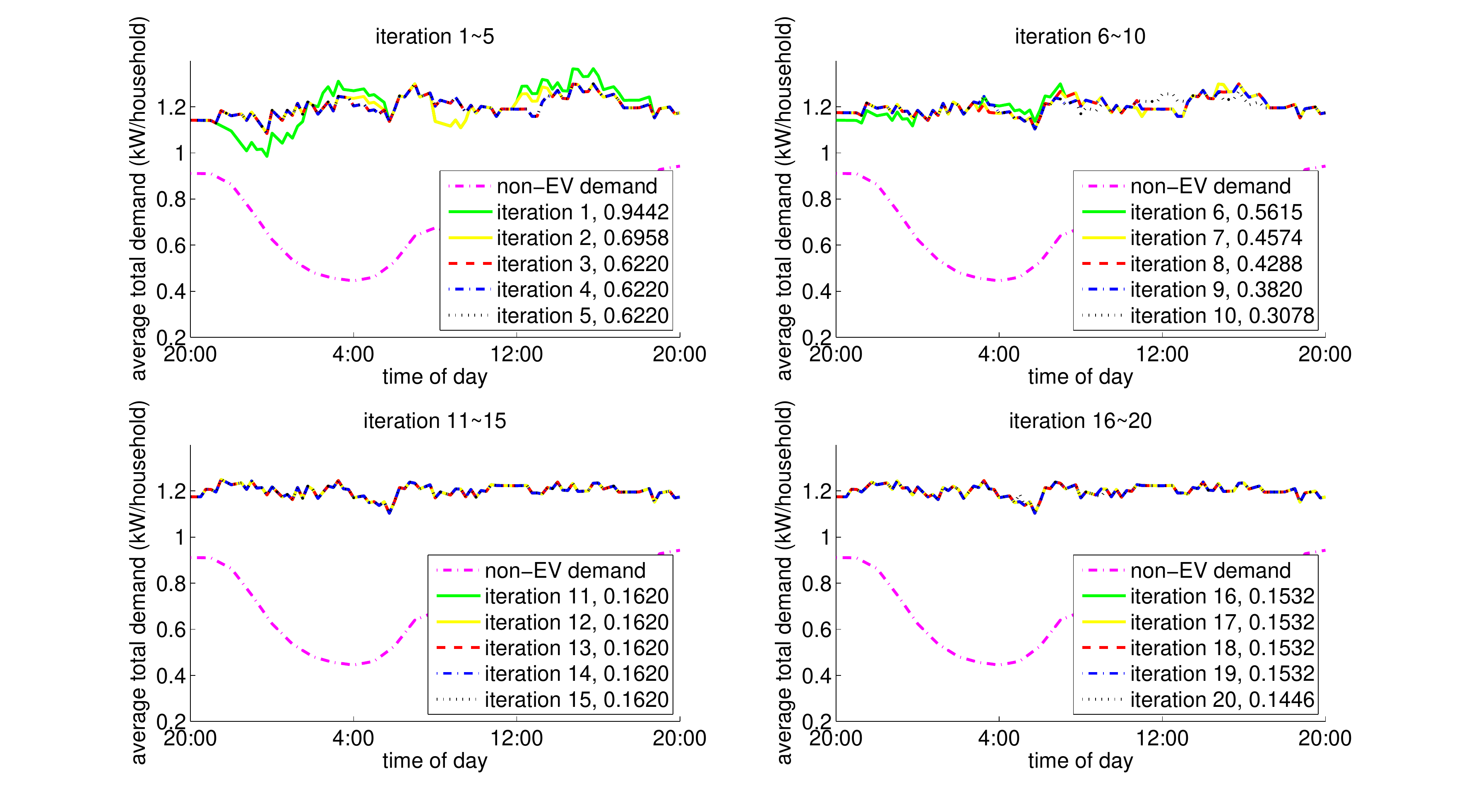}
      \caption{Average aggregate load profiles (per household) and escape probabilities in the first 20 iterations of Algorithm \ref{algorithm: nonconvex} in a 100\% EV penetration case. Escape probabilities for different iterations are shown in the legends.}
      \label{fig:convRate}
    \end{figure*}

To visualize the relationship between tentative charging profile updates and escape probabilities, we show the average aggregate load profiles (per household) and the escape probabilities in the first $20$ iterations of Algorithm \ref{algorithm: nonconvex} in a 100\% EV penetration case (one EV per household) in Figure \ref{fig:convRate}. It can be seen that the aggregate load profile only gets slightly updated from iteration $6$ to iteration $10$ where the escape probability ranges from $0.3$ to $0.6$; and that the aggregate load profile is updated only twice from iteration $11$ to iteration $20$ where the escape probability is below $0.2$. To summarize, the escape probability $P_\mathrm{escape}^{(k)}$ measures ``closeness'' to stationary charging profiles: the smaller $P_\mathrm{escape}^{(k)}$, the closer $x^{(k-1)}$ is to a stationary charging profile. In particular, if $P_\mathrm{escape}^{(k)}<0.5$, then $x^{(k-1)}$ is nearly stationary.

\begin{figure}[htbp!]
      \centering
      \includegraphics[scale=0.6]{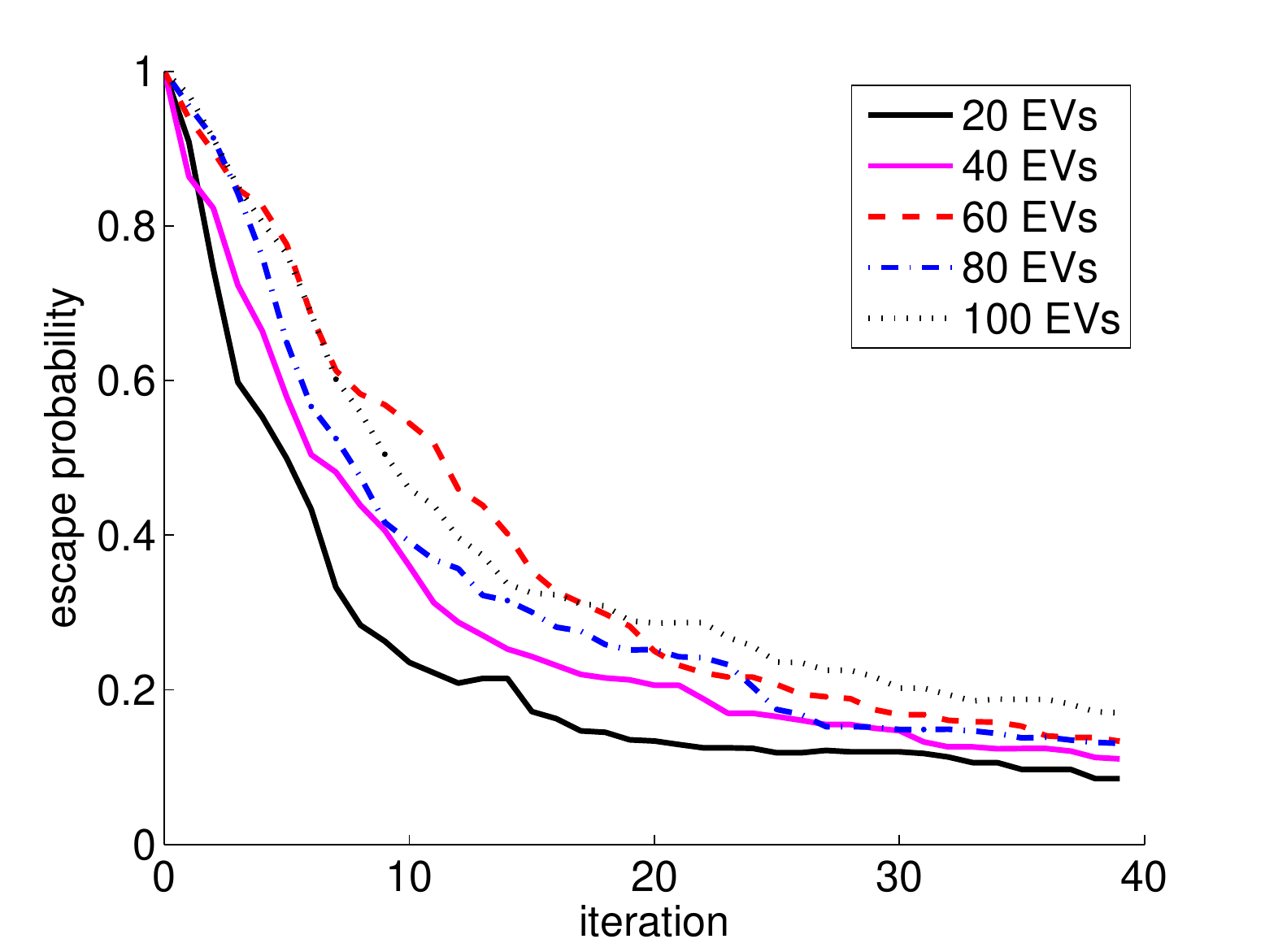}
      \caption{Average escape probability $P_\mathrm{escape}^{(k)}$ at different EV penetration levels.}
      \label{fig:escape}
    \end{figure}

Fig. \ref{fig:escape} shows the (average) escape probability $P_\mathrm{escape}^{(k)}$ (of $10$ simulations) in the first $40$ iterations of Algorithm \ref{algorithm: nonconvex} at different EV penetration levels. It can be seen that the average escape probability goes below $0.5$ within $20$ iterations for all penetration levels. Hence, we stop Algorithm \ref{algorithm: nonconvex} after $20$ iterations. We call $x^{(20)}$ the output charging profile of Algorithm \ref{algorithm: nonconvex} hereafter.

\subsection{Suboptimality ratio}\label{sec:simu:gap}
\begin{figure}[htbp!]
      \centering
      \includegraphics[scale=0.5]{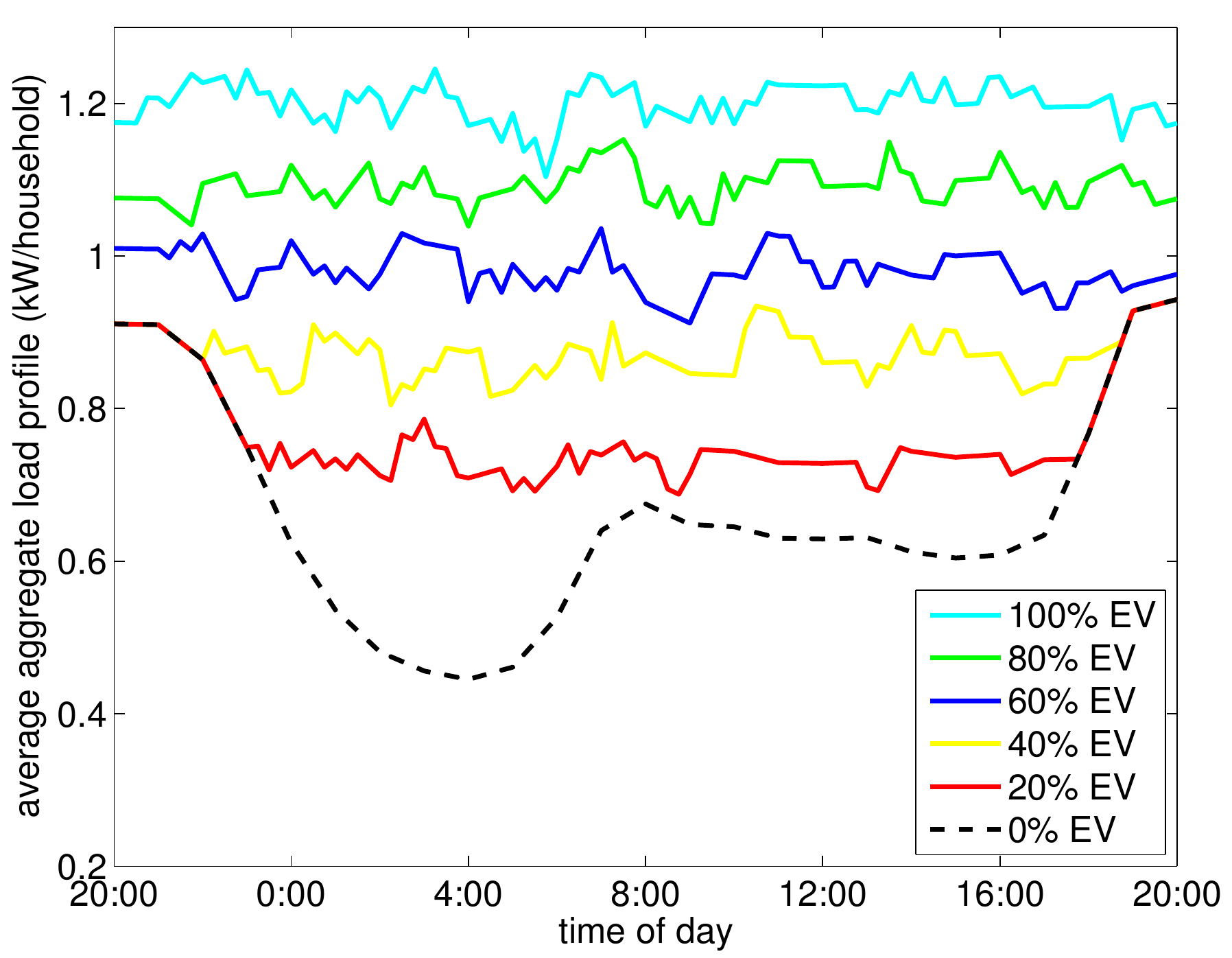}
      \caption{Average total demand per household in iteration $20$ of Algorithm DSC, with various number of EVs.}
      \label{fig:performance}
    \end{figure}

Fig. \ref{fig:performance} shows the average aggregate load profile in iteration $20$ of Algorithm \ref{algorithm: nonconvex} at different EV penetration levels. It can be seen that the average aggregate load profile is close to flat even with only 20\% EV penetration. Note that completely flat aggregate load profile is not achievable since the charging rate of an EV is either $0$ or $3.3$kW.

Now we theoretically quantify the suboptimality ratio $\subopt\left(x^{(20)}\right)$ [defined in \eqref{subopt}] of output charging profile $x^{(20)}$. Since tentative charging profile updates become negligible after $20$ iterations, we think of $x^{(20)}$ as a stationary charging profile and apply the suboptimality ratio upper bound derived in \eqref{subopt2}, i.e.,
	\begin{equation}
	\label{upper bound 1}
	\subopt\left(x^{(20)}\right)  \leq \frac{2\sum_{i=1}^n Y_i} {\|b+\sum_{i=1}^nx_i^*\|^2} := \mathrm{SubOptBound}
	\end{equation}
where $x^*$ denotes the optimal charging profile.

When the number $n$ of EVs is small, base load $b$ dominates aggregate EV load $\sum_{i=1}^nx_i^*$ and therefore
	\[ 2\sum_{i=1}^n Y_i ~=~ 2\sum_{i=1}^n\left\|x_i^*\right\|^2 ~\leq~ 2\left\|\sum_{i=1}^nx_i^*\right\|^2 ~\ll~ \left\|b\right\|^2 ~\leq~ \left\| b+\sum_{i=1}^nx_i^* \right\|^2. \]
Hence, $\mathrm{SubOptBound}$ is much smaller than 1 and the suboptimality ratio $\subopt\left(x^{(20)}\right)$ is small. When the number $n$ of EVs is big, base load $b$ is dominated by aggregate EV load and therefore
	\begin{equation}
	\label{upper bound 2}
	\mathrm{SubOptBound} ~\leq~ \frac{2\sum_{i=1}^nY_i}{\left\|\sum_{i=1}^nx_i^*\right\|^2}.
	\end{equation}
Noting that $2\sum_{i=1}^nY_i$ scales $O(n)$ as $n$ increases and that $\left\|\sum_{i=1}^nx_i^*\right\|^2$ scales $O(n^2)$ as $n$ increases, the bound $\mathrm{SubOptBound}$ scales $O(1/n)$ as $n$ increases. Hence, the suboptimality ratio $\mathrm{SubOpt}\left(x^{(20)}\right)$ remains small when $n$ is big.

\begin{figure}[htbp!]
      	\centering
      	\includegraphics[scale=0.4]{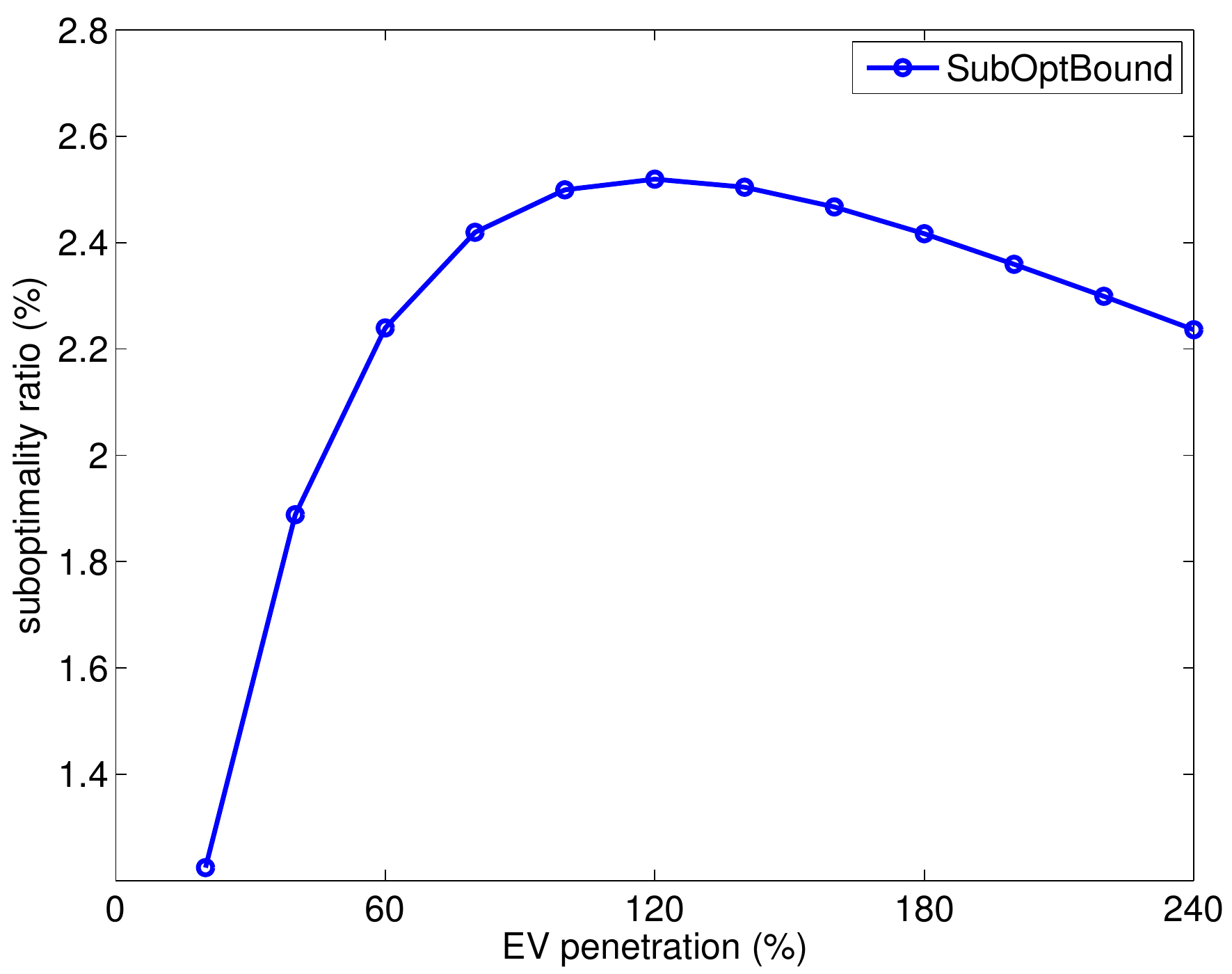}
      	\caption{Suboptimality upper bound $\mathrm{SubOptBound}$ for different EV penetration levels.}
      	\label{fig:gap}
	\end{figure}

The suboptimality upper bound $\mathrm{SubOptBound}$ for different penetration levels of EV is shown in Fig. \ref{fig:gap}. It can be seen that $\mathrm{SubOptBound}_1\leq2.6\%$ at all EV penetration levels.

\section{Conclusions}\label{sec:conclusion}
We have proposed a distributed EV charging scheduling algorithm to shape the aggregate load, that allows for a class of nonconvex constraints on admissible EV charging profiles. In particular, the algorithm applies to the case where each EV only has finitely many admissible charging profiles. The algorithm is iterative, and assumes the availability of a coordinator that can broadcast control signals to all EVs and receive feedbacks. In each iteration, the coordinator collects tentative charging profiles computed by the EVs in the previous iteration and broadcasts the normalized aggregate load profile. After receiving the broadcast signal, each EV updates its tentative charging profile by first computing a probability distribution over its admissible charging profiles, and then samples from the probability distribution to update its tentative charging profile.

The probability distribution an EV computes is designed such that the objective values in consecutive iterations of the algorithm decreases in expectation, and therefore the martingale theory ensures that the objective value converges almost surely as iterations continue. Besides, the algorithm shares the same information flow pattern as an algorithm that is designed to schedule EVs with convex constraints, and therefore both algorithms can run simultaneously in scenarios where there are both EVs with convex constraints and EVs with nonconvex constraints.

We have proved that the algorithm converges almost surely to a random stationary charging profile, and that the suboptimality ratio of the stationary charging profile has an upper bound that scales $O(1/n)$ as the number $n$ of EVs increases. Case studies confirm that the algorithm converges fast (can stop after 20 iterations), with suboptimality ratio below $2.6\%$ at all EV penetration levels.

\bibliographystyle{IEEEtran}
\bibliography{load_control}

\appendices
\section{Proof of Theorem \ref{thm: convergence convex}}
\label{app: convergence convex}
We apply the Lasalle Theorem to prove Theorem \ref{thm: convergence convex}. Note that $x^{(k)}$ depends continuously on $x^{(k-1)}$, the objective function
	\[ L(x)=\left\| b+\sum_{i=1}^nx_i \right\|^2 \]
is continuous in $x$, and the set $\prod_{i=1}^n\hX_i$ where $\{x^{(k)}:k\geq1\}$ take values in is compact. To prove Theorem \ref{thm: convergence convex}, i.e., $x^{(k)}\rightarrow \hO^*$ as $k\rightarrow\infty$, it suffices to prove (i) $L[x^{(k)}] \leq L[x^{(k-1)}]$ and (ii) $L[x^{(k)}] = L[x^{(k-1)}]$ implies $x^{(k-1)}\in\hO^*$ for $k\geq2$. The proofs of (i) and (ii) are presented below.

Fix an arbitrary $k\geq2$, and let $\Delta x_i:=x_i^{(k)}-x_i^{(k-1)}$ denote the update in $x_i$ in iteration $k$ of Algorithm \ref{algorithm: convex} for $i\in\hN$. Since $x_i^{(k)}$ minimizes \eqref{rate update} over $\hX_i$, it should have a smaller objective value than $x_i^{(k-1)}$, i.e.,
	\begin{align*}
	2c_i\left\langle g^{(k)}, x_i^{(k)} \right\rangle + \left\| x_i^{(k)}-x_i^{(k-1)}\right\|^2 \leq 2c_i\left\langle g^{(k)}, x_i^{(k-1)}\right\rangle 
	\qquad \Longrightarrow \qquad
	2c_i\left\langle g^{(k)}, \Delta x_i \right\rangle + \left\| \Delta x_i \right\|^2 \leq 0.
	\end{align*}
Furthermore, since \eqref{rate update} is strictly convex, the equality is attained if and only if $\Delta x_i=0$.

Abbreviate the aggregate service profile $d^{(k-1)}=b+\sum_{i=1}^nx_i^{(k-1)}$ in iteration $k-1$ by $d$ for convenience, then
    \begin{align*}
    L[x^{(k)}] - L[x^{(k-1)}]
    &=~ \left\| d+\sum_{i=1}^n\Delta x_i \right\|^2 - \left\| d \right\|^2 \\
    &=~ \left\| \sum_{i=1}^n\Delta x_i \right\|^2 + 2\sum_{i=1}^n \left\langle d, \Delta x_i\right\rangle \\
    &=~ \sum_{i=1}^n \left\| \Delta x_i \right\|^2 + \sum_{i\neq j} c_ic_j\left\langle \frac{\Delta x_i}{c_i}, \frac{\Delta x_j}{c_j} \right\rangle + 2\sum_{i=1}^n \left\langle d, \Delta x_i\right\rangle \\
    &\leq~ \sum_{i=1}^n \left\| \Delta x_i \right\|^2 + \frac{1}{2}\sum_{i\neq j} c_ic_j \left( \left\| \frac{\Delta x_i}{c_i} \right\|^2 + \left\| \frac{\Delta x_j}{c_j} \right\|^2 \right) + 2\sum_{i=1}^n \left\langle d, \Delta x_i\right\rangle \\
    &=~ \sum_{i=1}^n \left\| \Delta x_i \right\|^2 + \sum_{i\neq j} c_ic_j \left\| \frac{\Delta x_i}{c_i} \right\|^2 + 2\sum_{i=1}^n \left\langle d, \Delta x_i\right\rangle \\
    &=~ \sum_{i=1}^n \left\| \Delta x_i \right\|^2 + \sum_{i=1}^n \left\|\Delta x_i\right\|^2 \frac{\sum_{j\neq i} c_j}{c_i}     + 2\sum_{i=1}^n \left\langle d, \Delta x_i\right\rangle \\
    &=~ \sum_{i=1}^n \left\|\Delta x_i\right\|^2 \frac{\sum_{j=1}^n c_j}{c_i}     + 2\sum_{i=1}^n \left\langle d, \Delta x_i\right\rangle \\
    &=~ \sum_{i=1}^n\frac{\sum_{j=1}^nc_j}{c_i}\left( \left\| \Delta x_i\right\|^2
    +2c_i \left\langle g^{(k)}, \Delta x_i \right\rangle\right) ~\leq~ 0,
    \end{align*}
i.e., (i) holds.

Moreover, it follows from the inequality above that $L[x^{(k)}] = L[x^{(k-1)}]$ if and only if $\Delta x_i=0$ for $i\in\hN$. Therefore if $L[x^{(k)}] = L[x^{(k-1)}]$, then $x_i^{(k-1)}=x_i^{(k)}$ minimizes \eqref{rate update} and it follows from the first order optimality condition that
	\begin{align*}
	2c_i\left\langle g^{(k)},  ~x_i-x_i^{(k-1)}\right\rangle \geq 0, ~ i\in\hN,~x_i\in\hX_i
	\qquad \Longrightarrow \qquad
	2\left\langle d, ~x_i-x_i^{(k-1)}\right\rangle \geq 0, ~i\in\hN,~x_i\in\hX_i.
	\end{align*}
The latter is the first order optimality condition for $x^{(k-1)}$ to minimize the objective function $L(x)$. Hence, $L[x^{(k)}] = L[x^{(k-1)}]$ implies $x^{(k-1)}\in\hO^*$, i.e., (ii) holds.

We have proved that (i) and (ii) hold for an arbitrary $k\geq2$, and this completes the proof of Theorem \ref{thm: convergence convex}.

\section{Proof of Theorem \ref{thm: rate convex}}
\label{app: rate convex}
We prove \eqref{split} by mathematical induction on $k$.
\begin{itemize}
\item[i)] When $k=0$, $x_i^{(k)}=0$ for $i=0,\ldots,n$ and $x_{i'}^{(k)}=0$ for $i=1,\ldots,n$, therefore \eqref{split} holds.
\item[ii)] Assume that \eqref{split} holds for $k=K$ ($K\geq0$). We prove that \eqref{split} holds for $k=K+1$ as follows. It is straightforward to verify that $g^{(K+1)}$ is identical for two instances, and therefore
	\begin{eqnarray*}
	x_{i'}^{(K+1)}
	&=& \underset{x_{i'}\in\hX_{i'}} {\mathrm{argmin}}~ 2X_{i'}\left\langle g^{(K+1)}, x_{i'}\right\rangle + \left\| x_{i'}-x_{i'}^{(K)} \right\|^2\\
	&=& \underset{x_i\in\hX_i} {\mathrm{argmin}}~2X_i\left\langle g^{(K+1)}, x_i\right\rangle + \left\| x_i-x_i^{(K)} \right\|^2\\
	&=& x_i^{(K+1)}
	\end{eqnarray*}
for $i=2,\ldots,n$ and
	\begin{eqnarray*}
	x_{1'}^{(K+1)} &=& \underset{x_{1'}\in\hX_{1'}} {\mathrm{argmin}}~ 2X_{1'}\left\langle g^{(K+1)}, x_{1'}\right\rangle + \left\| x_{1'}-x_{1'}^{(K)} \right\|^2\\
	&=& \underset{x_{1'}\in2\hX_1} {\mathrm{argmin}}~ 4X_1\left\langle g^{(K+1)}, x_{1'}\right\rangle + \left\| x_{1'}-2x_1^{(K)} \right\|^2\\
	&=& 2\ \underset{2x_1\in2\hX_1} {\mathrm{argmin}}~ 4X_1\left\langle g^{(K+1)}, 2x_1\right\rangle + \left\| 2x_1-2x_1^{(K)} \right\|^2\\
	&=& 2\ \underset{x_1\in\hX_1} {\mathrm{argmin}}~ 2X_1\left\langle g^{(K+1)}, x_1\right\rangle + \left\| x_1-x_1^{(K)} \right\|^2\\
	&=& 2 x_1^{(K+1)}.
	\end{eqnarray*}
Similarly, $x_{1'}^{(K+1)} = 2X_0^{(K+1)}$. Hence, \eqref{split} holds for $k=K+1$.
\end{itemize}
According to (i) and (ii), \eqref{split} holds for $k\geq0$, which completes the proof of Theorem \ref{thm: rate convex}.

\section{Proof of Theorem \ref{thm: nash}}\label{app: nash}
To prove Theorem \ref{thm: nash}, we need to bridge a probability distribution $P_i$ over $x_i\in\hX_i$ with its expectation $\Ex_{P_i}[x_i]$. In particular, the following lemma is used when declaring a service profile $x$ is stationary for Algorithm \ref{algorithm: nonconvex}.

\begin{lemma}\label{lemma: delta}
Fix an $i\in\hN$ and assume that $\hX_i$ satisfies A4, i.e., every element in $\hX_i$ has the same $\ell_2$ norm $\sqrt{Y_i}$. Let $P_i$ be a probability distribution over $\hX_i$, then
	\[ \Ex_{P_i}[x_i]\in\hX_i
	\qquad \Longleftrightarrow \qquad
	P_i = \delta(x_i) \text{ for some }x_i\in\hX_i. \]
\end{lemma}

Lemma \ref{lemma: delta} implies that if $\Ex_{P_i^{(k)}}[x_i] = x_i^{(k-1)}$ for all $i$, then $P_i^{(k)}=\delta[x_i^{(k-1)}]$ for all $i$ and therefore the tentative charging profile $x^{(k-1)}$ is stationary for Algorithm \ref{algorithm: nonconvex}.

\begin{proof}
The ``$\Leftarrow$'' direction is straightforward, and the ``$\Rightarrow$'' direction can be proved as follows. One has
	\[ Y_i = \left\| \Ex_{P_i} [x_i] \right\|^2 \leq \Ex_{P_i} [\|x_i\|^2] = Y_i \]
and the equality is attained if and only if $P_i=\delta(x_i)$ for some $x_i\in\hX_i$. This completes the proof of Lemma \ref{lemma: delta}.
\end{proof}

\noindent{\it Proof of Theorem \ref{thm: nash}:}
According to Definition \ref{def: stationary}, a service profile $x^s=(x_1^s,\ldots,x_n^s)$ is stationary for Algorithm \ref{algorithm: nonconvex} if and only if for any $k\geq1$, $x^{(k-1)}=x^s$ implies $P_i^{(k)}=\delta(x_i^s)$ for $i\in\hN$.

($\Rightarrow$) Let $x^s$ be stationary for Algorithm \ref{algorithm: nonconvex}. Fix an arbitrary $k\geq1$ and assume $x^{(k-1)}=x^s$, then $P_i^{(k)}=\delta(x_i^s)$ for $i\in\hN$. Since $P_i^{(k)}$ minimizes \eqref{expectation} over $\Theta(\hX_i)$, $x_i^s=\Ex_{P_i^{(k)}}[x_i]$ minimizes
	\[ 2c_i\left\langle \frac{g^{(k)}\sum_{j=1}^nc_j-x_i^{(k-1)}}{\sum_{j\neq i}c_j}, x_i \right\rangle
	 + \left\| x_i - x_i^{(k-1)} \right\|^2 \]
over $x_i\in\conv(\hX_i)$ for $i\in\hN$. According to the first order optimality condition,
	\[ \left\langle 2c_i\frac{g^{(k)}\sum_{j=1}^nc_j-x_i^{(k-1)}}{\sum_{j\neq i}c_j} + 2\left(x_i^s-x_i^{(k-1)}\right), ~x_i-x_i^s\right\rangle \geq 0 \]
for all $x_i\in\conv(\hX_i)$ and all $i\in\hN$. Substitute $x^{(k-1)}=x^s$ and $g^{(k)}=\left(b+\sum_{j=1}^nx_j^{(k-1)}\right)/\sum_{j=1}^nc_j$ to obtain
	\[ \left\langle b+\sum_{j\neq i}x_j^s, ~x_i-x_i^s\right\rangle \geq 0 \]
for all $x_i\in\conv(\hX_i)$ and all $i\in\hN$. Then
	\begin{align*}
	\left\langle b+\sum_{j=1}^nx_j^s, ~x_i^s\right\rangle~
	&=~ \left\langle b+\sum_{j\neq i}x_j^s, ~x_i^s\right\rangle + Y_i \\
	&\leq~ \left\langle b+\sum_{j\neq i}x_j^s, ~x_i\right\rangle + Y_i \\
	&=~ \left\langle b+\sum_{j\neq i}x_j^s+x_i, ~x_i\right\rangle
	\end{align*}
for all $x_i\in\hX_i$ and all $i\in\hN$, i.e., $x^s$ is a Nash-equilibrium of the game described in Theorem \ref{thm: nash}.

($\Leftarrow$) Let $x^e$ be a Nash-equilibrium of the game described in Theorem \ref{thm: nash}, then
	\[ \left\langle b+\sum_{j=1}^nx_j^e, ~x_i^e\right\rangle \leq \left\langle b+\sum_{j\neq i}x_j^e+x_i, ~x_i\right\rangle \]
for all $x_i\in\hX_i$ and all $i\in\hN$. Since $\|x_i^e\|^2=Y_i=\|x_i\|^2$ for all $x_i\in\hX_i$, one has
	\begin{equation}
	\label{nash equilibrium}
	\left\langle b+\sum_{j\neq i}x_j^e, ~x_i^e \right\rangle \leq \left\langle b+\sum_{j\neq i}x_j^e, ~x_i \right\rangle
	\end{equation}	
for all $x_i\in\hX_i$ and all $i\in\hN$. It follows that \eqref{nash equilibrium} also holds for all $x_i\in\conv(\hX_i)$ and all $i\in\hN$. Then,
	\[ \left\langle b+\sum_{j\neq i}x_j^e, ~x_i-x_i^e\right\rangle \geq 0
	\quad\Longrightarrow\quad
	\left\langle 2c_i\frac{b+\sum_{j\neq i}x_j^e}{\sum_{j\neq i}c_j}  +   2\left(x_i^e-x_i^e\right), ~x_i-x_i^e\right\rangle \geq 0 \]	
for all $x_i\in\conv(\hX_i)$ and all $i\in\hN$. If $x^{(k-1)}=x^e$ for some $k\geq1$, then
	\[ \left\langle 2c_i\frac{g^{(k)}\sum_{j=1}^nc_j-x_i^{(k-1)}}{\sum_{j\neq i}c_j}   +   2\left(x_i^e-x_i^{(k-1)}\right), ~x_i-x_i^e\right\rangle \geq 0 \]
for all $x_i\in\conv(\hX_i)$ and all $i\in\hN$. According to the first order optimality condition, $x_i^e$ minimizes
	\begin{equation}
	\label{obj: expectation}
	L(x_i) := 2c_i\left\langle \frac{g^{(k)}\sum_{j=1}^nc_j-x_i^{(k-1)}}{\sum_{j\neq i}c_j}, x_i \right\rangle
	 + \left\| x_i - x_i^{(k-1)} \right\|^2
	 \end{equation}
over $x_i\in\conv(\hX_i)$ for $i\in\hN$. Since \eqref{obj: expectation} is continuous in $x_i$, $x_i^e$ also minimizes \eqref{obj: expectation} over $\overline{\conv(\hX_i)}$ for $i\in\hN$.

For each $i\in\hN$ and each probability distribution $P_i$ over $\hX_i$, let
	\[ L(P_i) := L(\Ex_{P_i}[x_i]) \]
denote the value of \eqref{obj: expectation} evaluated at the expectation $\Ex_{P_i}[x_i]$ of $P_i$, then $L(P_i)$ coincides with the objective in \eqref{expectation}. On one hand, since $\hX_i$ is compact, the expectation $\Ex_{P_i}[x_i]\in\overline{\conv(\hX_i)}$ and therefore
	\[ L(P_i) \geq \min_{x_i\in\overline{\conv(\hX_i)}} L(x_i) = L(x_i^e) \]
for all $P_i\in\Theta(\hX_i)$. On the other hand, the distribution $\delta(x_i^e)$ achieves
	\[ L[\delta(x_i^e)] = L(x_i^e). \]
Hence, the probability distribution $\delta(x_i^e)$ minimizes \eqref{expectation} over $\Theta(\hX_i)$.

Moreover, $\delta(x_i^e)$ is the unique minimizer of $L(P_i)$ over $P_i\in\Theta(\hX_i)$. Otherwise, there exists $P_i'$ other than $\delta(x_i^e)$ that minimizes $L(\cdot)$ over $\Theta(\hX_i)$. It follows that $L(P_i')=L(\Ex_{P_i'}[x_i])=L(x_i^e)$, i.e., both $x_i^e$ and $\Ex_{P_i'}[x_i]$ minimize $L(x_i)$ over $x_i\in\overline{\conv(\hX_i)}$. Since $L(x_i)$ is strictly convex in $x_i$, one must have
	\[ \Ex_{P_i'}[x_i] = x_i^e \in \hX_i. \]
It follows from Lemma \ref{lemma: delta} that $P_i'=\delta(x_i)$ for some $x_i\in\hX_i$ and therefore $P_i'=\delta(x_i^e)$. This contradicts our assumption that $P_i'\neq\delta(x_i^e)$. Hence, $\delta(x_i^e)$ is the unique minimizer of $L(P_i)$ over $P_i\in\Theta(\hX_i)$.

Hence, $P_i^{(k)} = \delta(x_i^e)$ for $i\in\hN$ and therefore $x^e$ is stationary.

Combining ($\Rightarrow$) and ($\Leftarrow$) completes the proof of Theorem \ref{thm: nash}.
$\hfill\blacksquare$

\section{Proof of Theorem \ref{thm: suboptimality}}\label{app: suboptimality}
The gap $\|b+\sum_{i=1}^nx_i^*\|^2-\|b+\sum_{i=1}^nx_i^s\|^2$ can be bounded below as follows.
    \begin{eqnarray*}
    \left\| b+\sum_{i=1}^n x_i^* \right\|^2 - \left\| b+\sum_{i=1}^n x_i^s \right\|^2
    &=& 2 \left\langle b+\sum_{i=1}^n x_i^s, ~\sum_{i=1}^n (x_i^*-x_i^s) \right\rangle + \left\| \sum_{i=1}^nx_i^* - \sum_{i=1}^nx_i \right\|^2 \\
    &\geq& 2\sum_{i=1}^n \left\langle b+\sum_{j=1}^n x_j^s, ~x_i^*-x_i^s \right\rangle\\
    &=& 2\sum_{i=1}^n \left\langle b+\sum_{j\neq i} x_j^s, ~x_i^*-x_i^s \right\rangle
    + 2\sum_{i=1}^n \left\langle x_i^s, ~x_i^*-x_i^s \right\rangle \\
    &\geq& 2\sum_{i=1}^n\left\langle x_i^s, x_i^*-x_i^s \right\rangle 			\qquad\qquad \text{[by \eqref{nash equilibrium}]}\\
    &=& 2\sum_{i=1}^n\left\langle x_i^s, x_i^* \right\rangle -2\sum_{i=1}^nY_i.
    \end{eqnarray*}
Since $|\left\langle x_i^s, x_i^* \right\rangle| \leq \|x_i^s\|\cdot\|x_i^*\| = \sqrt{Y_i}\cdot\sqrt{Y_i}=Y_i$, one has
	$ \left\langle x_i^s, x_i^*\right\rangle \geq - Y_i $
for $i\in\hN$ and therefore
	\[ \left\| b+\sum_{i=1}^n x_i^* \right\|^2 - \left\| b+\sum_{i=1}^n x_i^s \right\|^2
	\geq 2\sum_{i=1}^n\left\langle x_i^s, x_i^* \right\rangle -2\sum_{i=1}^nY_i
	\geq -4\sum_{i=1}^nY_i. \]
Further, if $\hX_i$ is composed of nonnegative elements, then $\left\langle x_i^s, x_i^* \right\rangle\geq0$ for $i\in\hN$ and therefore
	\[ \left\| b+\sum_{i=1}^n x_i^* \right\|^2 - \left\| b+\sum_{i=1}^n x_i^s \right\|^2
	\geq 2\sum_{i=1}^n\left\langle x_i^s, x_i^* \right\rangle -2\sum_{i=1}^nY_i
	\geq -2\sum_{i=1}^nY_i. \]
This completes the proof of Theorem \ref{thm: suboptimality}.

\section{Proof of Proposition \ref{proposition: supermartingale}}\label{app: supermartingale}
Fix an arbitrary $k\geq2$ and let $\Delta x_i:=x_i^{(k)}-x_i^{(k-1)}$ denote the update of $x_i$ in iteration $k$ of Algorithm \ref{algorithm: nonconvex} for $i\in\hN$. Note that the process $\{x^{(k)} : k\geq0 \}$ is Markov. Therefore
	\[ \Ex\left[ f[x^{(k)}] \mid x^{(k-1)},x^{(k-2)},\ldots,x^{(0)} \right] = \Ex\left[ f[x^{(k)}] \mid x^{(k-1)} \right] \]
for any function $f$. Abbreviate $\Ex\left[\Delta x_i \mid x^{(k-1)}\right]$ by $\Ex\Delta x_i$ and abbreviate $\Ex\left[ x_i^{(k)}\mid x^{(k-1)}\right]$ by $\Ex x_i^{(k)}$ for $i\in\hN$ without confusion. Also abbreviate the aggregate service profile $d^{(k-1)}=b+\sum_{i=1}^nx_i^{(k-1)}$ in iteration $k-1$ of Algorithm \ref{algorithm: nonconvex} by $d$ and abbreviate $g^{(k)}=d/\sum_{i=1}^nc_i$ by $g$, then
    \begin{align*}
    \Ex\left[ L_k \mid x^{(k-1)} \right] -L_{k-1}
    ~&=~ \Ex\left[ \left\|d+\sum_{i=1}^n\Delta x_i \right\|^2 \mid x^{(k-1)} \right] - \left\|d \right\|^2\\
    &=~ \Ex \left[ \left\|\sum_{i=1}^n\Delta x_i \right\|^2 \mid x^{(k-1)} \right]
    + 2\sum_{i=1}^n     \left\langle d, \Ex\Delta x_i    \right\rangle.
    \end{align*}
The first term can be simplified as
    \begin{align*}
    \Ex \left[ \left\|\sum_{i=1}^n\Delta x_i \right\|^2 \mid x^{(k-1)} \right]
    ~&=~ \Ex \left[        \sum_{i=1}^n\left\|\Delta x_i \right\|^2    +      \sum_{i\neq j}\left\langle \Delta x_i, \Delta x_j \right\rangle        \mid x^{(k-1)} \right]\\
    &=~ \sum_{i=1}^n\Ex \left[        \left\|\Delta x_i \right\|^2     \mid x^{(k-1)}\right]
    + \sum_{i\neq j}           \left\langle \Ex \Delta x_i , \Ex \Delta x_j \right\rangle        \\
    \end{align*}
where
    \begin{align*}
    \Ex \left[       \left\|\Delta x_i \right\|^2     \mid x^{(k-1)} \right]~
    &=~ \Ex \left[        \left\| x_i^{(k)}-x_i^{(k-1)} \right\|^2     \mid x^{(k-1)} \right] \\
    &=~ \Ex \left[        \left\| x_i^{(k)} \right\|^2     \mid x^{(k-1)} \right]
    +\Ex \left[        \left\| x_i^{(k-1)} \right\|^2     \mid x^{(k-1)}\right]
    -2 \Ex \left[        \left\langle x_i^{(k)}, x_i^{(k-1)} \right\rangle     \mid x^{(k-1)}\right]\\
    &=~ 2Y_i      -2 \left\langle \Ex x_i^{(k)}, x_i^{(k-1)} \right\rangle    ~=~ -2 \left\langle \Ex\Delta x_i, x_i^{(k-1)} \right\rangle\\
    \end{align*}
for $i=1,2,\ldots,n$ and
    \begin{eqnarray*}
    \sum_{i\neq j} \left\langle \Ex \Delta x_i , \Ex \Delta x_j \right\rangle
    &=& \sum_{i\neq j} c_ic_j \left\langle \frac{\Ex \Delta x_i}{c_i} , \frac{\Ex \Delta x_j}{c_j} \right\rangle \\
    &\leq& \sum_{i\neq j} \frac{c_ic_j}{2}       \left( \frac{\left\|\Ex \Delta x_i \right\|^2}{c_i^2} +  \frac{\left\|\Ex \Delta x_j\right\|^2}{c_j^2}                         \right)        \\
    &=& \sum_{i\neq j}     c_ic_j       \frac{\left\|\Ex \Delta x_i \right\|^2}{c_i^2}\\
    &=& \sum_{i=1}^n  \frac{\sum_{j\neq i}c_j}{c_i}       \left\|\Ex \Delta x_i\right\|^2.
    \end{eqnarray*}
Therefore
    \begin{align}
    \Ex\left[ L_k \mid x^{(k-1)} \right]-L_{k-1}~
    &=~ \sum_{i=1}^n\left[ -2 \left\langle \Ex\Delta x_i, x_i^{(k-1)} \right\rangle
    +\frac{\sum_{j\neq i}c_j}{c_i}       \left\|\Ex \Delta x_i\right\|^2  + 2\left\langle d, \Ex\Delta x_i    \right\rangle\right] \nonumber\\
    &=~ \sum_{i=1}^n\left[ \frac{\sum_{j\neq i}c_j}{c_i}       \left\|\Ex \Delta x_i\right\|^2  + 2\left\langle d-x_i^{(k-1)}, \Ex\Delta x_i    \right\rangle\right] \nonumber\\
    &=~ \sum_{i=1}^n \frac{\sum_{j\neq i}c_j}{c_i}     \left[        \left\|\Ex \Delta x_i\right\|^2  + 2c_i\left\langle \frac{d-x_i^{(k-1)}}{\sum_{j\neq i} c_j}, ~\Ex\Delta x_i    \right\rangle\right] \nonumber\\
    &=~ \sum_{i=1}^n \frac{\sum_{j\neq i}c_j}{c_i}     \left[        \left\|\Ex \Delta x_i\right\|^2  + 2c_i\left\langle \frac{g\sum_{j=1}^n c_j -x_i^{(k-1)}}{\sum_{j\neq i} c_j}, ~\Ex\Delta x_i    \right\rangle\right] \nonumber\\
    &\leq~0 \label{condition}
    \end{align}
The last inequality is due to the fact that $P_i^{(k)}$ minimizes \eqref{expectation} and therefore should have a smaller objective value than $\delta(x_i^{k-1})$. This completes the proof of Proposition \ref{proposition: supermartingale}.

\section{Proof of Theorem \ref{thm: convergence}}\label{app: convergence}
The key of proving Theorem \ref{thm: convergence} is the following lemma, which generalizes the Lassale Theorem from deterministic processes to ``hidden'' lower bounded supermartingales.

\begin{lemma}
\label{lemma: lassale martingale}
Let $\{X_k:k\geq1\}$ be a time-invariant Markov process that takes values on a compact set $\hX$, i.e., $X_k$ is a map from the sample space to $\hX$ for $k\geq1$. Let $L:\hX\mapsto\mathbb{R}$ be bounded, i.e., there exists $M\in\mathbb{R}$ such that $|L(x)| \leq M$ for $x\in\hX$. Assume that $\{L(X_k):k\geq1\}$ is a supermartingale with respect to $\{X_k:k\geq1\}$, i.e.,
	\[ \Ex\left[L(X_k) \mid X_{k-1}\right] \leq L(X_{k-1}), \qquad k\geq2. \]
Also assume that the function $f:\hX\mapsto\mathbb{R}_-$ defined by
	\[ f(x) := \Ex\left[L(X_k) \mid X_{k-1}=x\right] - L(x) \leq 0 \]
is continuous and let
	\[ \hS:=\{x\in\hX \mid f(x)=0\} \]
denote the set of $x\in\hX$ where $\Ex[L(X_k)\mid X_{k-1}=x] = L(x)$. Then
	\[ \dist(X_k,\hS) \overset{\text{a.s.}}{\longrightarrow} 0, \qquad k\rightarrow \infty. \]
\end{lemma}

Lemma \ref{lemma: lassale martingale} implies that if we can find a function $L$ such that the stochastic process $L(X_k)$ decreases in expectation, then the stochastic process $X_k$ almost surely converges to certain set $\hS$, under some technical conditions. Note that if the process $\{X_k:k\geq1\}$ is deterministic, then Lemma \ref{lemma: lassale martingale} reduces to the Lassale Theorem. The stochastic process $\{X_k:k\geq1\}$ is called a hidden lower bounded supermartingale since its function $L(X_k)$ is a lower bounded supermartingale. The proof of Lemma \ref{lemma: lassale martingale} is similar to the proof of the Martingale Convergence Theorem \cite[Chapter 12.3]{grimmett2001probability} and presented below.

\begin{proof}
The statement can be transformed as follows.
    \begin{align*}
    \pr\left\{\lim_{k\rightarrow\infty}\dist(X_k,\hS)=0\right\}=1~
    &\Longleftrightarrow~ \pr \left\{\underset{n\geq1}{\cap}\left\{\underset{k\rightarrow\infty}{\overline{\lim}}\dist(X_k,\hS) <\frac{1}{n}\right\}\right\}=1\\
    &\Longleftrightarrow~ \lim_{n\rightarrow\infty}\pr \left\{\underset{k\rightarrow\infty}{\overline{\lim}}\dist(X_k,\hS) <\frac{1}{n}\right\}=1\\
    &\Longleftrightarrow~ \pr \left\{\underset{k\rightarrow\infty}{\overline{\lim}}\dist(X_k,\hS) <\frac{1}{n}\right\}=1,~\forall n\geq1\\
    &\Longleftrightarrow~ \pr \left\{\dist(X_k,\hS) \geq\frac{1}{n}\text{ infinitely often} \right\}=0,~\forall n\geq1\\
    &\Longleftrightarrow~ \pr \left\{\#\left\{k \mid \dist(X_k,\hS) \geq\frac{1}{n}\right\}=\infty \right\}=0,~\forall n\geq1\\
    &\Longleftarrow~ \Ex \left[\#\left\{k \mid \dist(X_k,\hS) \geq\frac{1}{n}\right\} \right]<\infty,~\forall n\geq1.
    \end{align*}
In the rest of the proof, we show that $\Ex \left[\#\left\{k \mid \dist(X_k,\hS) \geq 1/n \right\} \right]<\infty$ for $n\geq1$.
    
Fix an arbitrary $n\geq1$ and let $T:=\#\left\{k\mid\dist(X_k,\hS) \geq1/n\right\}$ denote the number of iterations that $X_k$ is at least $1/n$ distance from $\hS$, then we want to prove $\Ex[T]<\infty$. Define an auxiliary process $\{\mathbbm{1}_k:k\geq1\}$ by
    \[ \mathbbm{1}_k:=\begin{cases}
    1 & \text{if }\dist(X_k,\hS)\geq 1/n,\\
    0 & \text{otherwise}.
    \end{cases} \]
Since the set $\mathbb{A}_n:=\{x\in\hX \mid \dist(x,\hS)\geq 1/n \}$ is compact, $f$ is continuous, and $f<0$ on $\mathbb{A}_n$, there exists $\epsilon>0$ such that
	\[ f(x) \leq -\epsilon, \qquad x\in\mathbb{A}_n. \]
Hence,
	\begin{align*}
	\mathbbm{1}_{k-1}=1 \quad
	&\Rightarrow\quad \dist(X_{k-1},\hS)\geq1/n \\
	&\Rightarrow\quad X_{k-1} \in \mathbb{A}_n \\
	&\Rightarrow\quad f(X_{k-1}) \leq -\epsilon
	\end{align*}
for $k\geq2$ and therefore
    \begin{align*}
    \Ex[L(X_K)-L(X_1)]~
    &=~ \sum_{k=2}^K \Ex[L(X_k)-L(X_{k-1})] ~=~ \sum_{k=2}^K \Ex\left[f(X_{k-1})\right]\\
    &\leq~ \sum_{k=2}^K \Ex\left[\mathbbm{1}_{k-1}f(X_{k-1})\right] ~\leq~ \sum_{k=2}^K \Ex\left[\mathbbm{1}_{k-1}\cdot(-\epsilon)\right]\\
    &=~ -\epsilon \sum_{k=1}^{K-1}  \Ex[\mathbbm{1}_k]
    \end{align*}
for $K\geq1$.
Finally,
    \[ \Ex [T]
    = \Ex \left[\sum_{k\geq1}\mathbbm{1}_k\right]
    =\lim_{K\rightarrow\infty}\sum_{k=1}^{K-1}\Ex \left[\mathbbm{1}_k\right]
    \leq\lim_{K\rightarrow\infty}\frac{\Ex[L(X_1)]-\Ex[L(X_K)]}{\epsilon}
    \leq\lim_{K\rightarrow\infty}\frac{2M}{\epsilon}<\infty. \]
This completes the proof of Lemma \ref{lemma: lassale martingale}.
\end{proof}

Now we apply Lemma \ref{lemma: lassale martingale} to the proof of Theorem \ref{thm: convergence}. In particular, let $X_k = x^{(k)}$ for $k\geq1$, and let $L$ be the objective function, then $L[x^{(k)}]=L_k$ is a supermartingale. In order to apply Lemma \ref{lemma: lassale martingale}, we are left to prove that
	\[ \Ex\left[ L_{k} \mid x^{(k-1)} \right] = L_{k-1} ~\Longleftrightarrow~ x^{(k-1)} \text{ is stationary.} \]
Fix an arbitrary $k\geq2$, then
    \begin{eqnarray*}
    \Ex\left[ L_{k} \mid x^{(k-1)} \right] = L_{k-1}
    &\Longleftrightarrow& \Ex[ x_i^{(k)} \mid x_i^{(k-1)} ] = x_i^{(k-1)} \text{ for }i\in\hN  \qquad\qquad \text{[by \eqref{condition}]} \\
    &\Longleftrightarrow& \Ex_{P_i^{(k)}}[x_i]=x_i^{(k-1)} \text{ for }i\in\hN \\
    &\Longleftrightarrow& P_i^{(k)}=\delta\left[x_i^{(k-1)}\right] \text{ for }i\in\hN 			\qquad\qquad \text{[by Lemma \ref{lemma: delta}]} \\
    &\Longleftrightarrow& x_i^{(k)}=x_i^{(k-1)} \text{ for }i\in\hN \\
    &\Longleftrightarrow& x^{(k-1)} \text{ is stationary}.
    \end{eqnarray*}
This completes the proof of Theorem \ref{thm: convergence}.

\section{Proof of Corollary \ref{cor: convergence to points}}\label{app: convergence to points}
Since the set $\hX:=\prod_{i=1}^n\hX_i$ has finitely many elements, the minimum distance
	\[ d_{\min} := \min_{x^1,x^2\in\hX} \| x^1 - x^2 \| \]
among different elements in $\hX$ is strictly positive. It follows from Theorem \ref{thm: convergence} that $\dist(x^{(k)},\hS)<d_{\min}$ almost surely as $k\rightarrow\infty$.

Pick a random sample path of the process $\{x^{(k)}: k\geq1\}$. With probability 1, there exists $k_0\geq1$ such that $\dist(x^{(k_0)},\hS)<d_{\min}$. It must be that $x^{(k_0)}\in\hS$ since otherwise the distance is at least $d_{\min}$. Then, since elements in $\hS$ are stationary service profiles, one must have $x^{(k_0)}=x^{(k_0+1)}=\cdots$ and therefore $\lim_{k\rightarrow\infty}x^{(k)}=x^{(k_0)}\in\hS$ on this sample path.

The sample pathes $\{x^{(k)}:k\geq1\}$ may converge (with probability 1) to different $x^{(k_0)}\in\hS$. Denote the random service profile that the process $\{x^{(k)}:k\geq1\}$ converges to by $x^\infty$, then
	\[ x^{(k)} \overset{\text{a.s.}}{\longrightarrow} x^\infty \text{ as }k\rightarrow\infty. \]

\end{document}